\documentclass{paper}

\usepackage{hyperref}
\usepackage{graphicx}
\usepackage{amsmath}
\usepackage{xcolor}
\usepackage{amssymb}
\usepackage{amsfonts}
\usepackage{marginnote}
\usepackage{fancybox}
\usepackage{color}
\usepackage{bm}
\usepackage{cancel}
\usepackage{wrapfig}
\usepackage{lipsum}
\usepackage{epstopdf}
\usepackage{stmaryrd}
\usepackage{mathrsfs}
\usepackage{psfrag}
\usepackage{caption}
\usepackage{mathtools}
\usepackage{amsopn}
\usepackage{stmaryrd}
\usepackage{placeins}
\usepackage[normalem]{ulem}

\newtheorem{remark}{Remark}[section]

\newcommand{\mbRn}{{\mathbb{R}^n}}
\newcommand{\T}{\mathscr{T}}
\newcommand{\V}{\mathbb{V}}

\newcommand{\omg}{{\Omega}}
\newcommand{\omgc}{{\Omega^c}}

\newtheorem{theorem}{Theorem}
\newtheorem{definition}{Definition}
\newtheorem{lemma}{Lemma}
\newtheorem{corollary}{Corollary}
\newtheorem{proposition}{Proposition}

\DeclareMathOperator*{\diam}{diam}
\renewcommand{\tilde}[1]{\widetilde{#1}}

\title{A priori error estimates for the optimal control of the integral fractional Laplacian\thanks{Submitted to the editors \today. }}
\author{Marta D'Elia\thanks{Center for Computing Research, Sandia National Laboratories, Albuquerque, NM 87123, USA ({\tt mdelia@sandia.gov}).}
\and
Christian Glusa\thanks{Center for Computing Research, Sandia National Laboratories, Albuquerque, NM 87123, USA ({\tt caglusa@sandia.gov}).}
\and
Enrique Ot\'arola\thanks{Departamento de Matem\'atica, Universidad T\'ecnica Federico Santa Mar\'ia, Valpara\'iso, Chile ({\tt enrique.otarola@usm.cl}).}
}

\begin{document}
\maketitle

\begin{abstract}
We design and analyze solution techniques for a linear-quadratic optimal control problem involving the integral fractional Laplacian. We derive existence and uniqueness results, first order optimality conditions, and regularity estimates for the optimal variables. We propose two strategies to discretize the fractional optimal control problem: a semidiscrete approach where the control is not discretized -- the so-called variational discretization approach -- and a fully discrete approach where the control variable is discretized with piecewise constant functions. Both schemes rely on the discretization of the state equation with the finite element space of continuous piecewise polynomials of degree one. We derive a priori error estimates for both solution techniques. We illustrate the theory with two-dimensional numerical tests.
\end{abstract}


\medskip\noindent\textbf{\textsf{AMS}} 35R11, 49J20, 49M25, 65K10, 65N15, 65N30.

\section{Introduction}
Nonlocal models have recently become of great interest to the applied sciences and engineering.
This is mainly due to the fact that operators featuring nonlocal interactions better describe many processes, for instance, anomalous diffusion phenomena, for which where classical integer order differential operators fail to provide an accurate description.
More specifically, they arise in applications such as stochastic jump processes \cite{XXXX,Valdinoci2009_FromLongJumpRandom}, material science (e.g. subsurface flow where nonlocal porous media models accurately describe the physical process) \cite{BensonWheatcraftEtAl2000_ApplicationFractionalAdvectionDispersionEquation,Biler2015,Silling2000_ReformulationElasticityTheoryDiscontinuities}, image processing \cite{GilboaOsher2008_NonlocalOperatorsWithApplications,MR2578033}, finance \cite{MR2064019,Tankov2003_FinancialModellingWithJumpProcesses}, fluids \cite{doi:10.1063/1.2208452,doi:10.1063/1.2819487,MccayNarasimhan1981_TheoryNonlocalElectromagneticFluids}, population dynamics \cite{Clark_2001} and cardiology \cite{Cusimano_2015}.

Fractional operators are a particular class of nonlocal operators.
When solving fractional partial differential equations (PDEs), even linear ones, several modeling and computational challenges arise.
As an example, the computational cost required by the solution of a linear fractional PDE can be prohibitively expensive, especially in two- or three-dimensional domains.
This is due to the fact that, contrary to the case of local PDEs, points in a domain interact with every other point in the space, due to the nonlocal nature of the operator that allows for infinite range interactions. This clearly creates computational challenges as the discretized problems are hard to assemble and solve.

Furthermore, it is often the case that the mathematical model is not exact, e.g. source terms or coefficients may be unknown or subject to uncertainty.
However, in the case when limited data or a priori information is available, one can resort to the solution of a control or inverse problem to recover the unknown parameters and define a more accurate, data-driven, mathematical model.

Among available data we may have sparse and/or noisy measurements of the state of the system or of an output of interest that we would like to match.
In this work we address the problem of finding an input function (e.g. a distributed source term) such that the corresponding solution is as close a possible to a target state; for now, we do not consider any uncertainty in the data.
We propose to solve an optimal control problem where the cost functional quantifies the misfit between the target and the predicted output of interest, the constraint is the fractional differential equation, and the control is a distributed source term.

PDE-constrained optimization problems involving fractional and nonlocal equations are not new in the literature; we mention, e.g., the works by Antil and Ot\'arola \cite{MR3429730}, Ot\'arola \cite{MR3702421}, and D'Elia and Gunzburger \cite{MR3158780,Delia2016}. In \cite{MR3429730}, the authors consider a linear-quadratic optimal control problem for the spectral definition of the fractional Laplacian; control constraints are also considered. The authors also propose and study solution techniques to approximate the underlying solution. In \cite{MR3158780} the authors consider an optimal control problem for a general nonlocal diffusion operator with finite range interactions. In the current work, with a similar formulation, we consider a linear-quadratic optimal control problem involving the integral definition of the fractional Laplace operator, which we simply refer to as the {\it integral fractional Laplacian}; in this case, as previously mentioned, the interactions can be infinite. It is important to note that the integral and spectral definitions of the fractional Laplace operator \emph{do not coincide}. In fact, as shown in \cite{MR3246044}, their difference is positive and positivity preserving. This, in particular, implies that the boundary behavior of the solutions to basic problems involving the aforementioned definitions are different \cite{MR3489634,MR1204855,MR3276603}.

In this work, we design and analyze efficient solution techniques for a linear-quadratic optimal control problem involving the integral fractional Laplacian. To make matters precise, for $n \geq 1$, we let $\Omega \subset \mathbb{R}^n$ be an open bounded domain with Lipschitz boundary $\partial \Omega$. Given $s \in (0,1)$ and a desired state $u_d: \Omega \rightarrow \mathbb{R}$, we define the cost functional
\begin{equation}
 \label{eq:cost_functional}
 J(u,z) := \frac{1}{2} \| u- u_d \|_{L^2(\Omega)}^2 + \frac{\alpha}{2} \| z \|_{L^2(\Omega)}^2,
\end{equation}
where $\alpha > 0$ denotes the so-called regularization parameter. Let $f :\Omega \rightarrow \mathbb{R}$ be a fixed function. We consider the following optimal control problem: Find
\begin{equation}
\label{eq:min}
 \min  J(u,z)
\end{equation}
subject to the \emph{fractional state equation}
\begin{equation}
\label{eq:state_equation}
(-\Delta)^s u = f + z \textrm{ in } \Omega, \qquad u = 0 \textrm{ in } \Omega^c,
\end{equation}
with $s\in(0,1)$ and $\omgc=\mbRn\setminus \overline\omg$, and the \emph{control constraints}
\begin{equation}
 \label{eq:control_constraints}
 a \leq z(x) \leq b \quad \textrm{a.e.} \quad x \in \Omega;
\end{equation}
the control bounds $a,b \in \mathbb{R}$ are such that $a < b$.

For functions defined over the whole space $\mathbb{R}^n$, the integral fractional Laplacian $(-\Delta)^s$ can be naturally defined via the Fourier transform as follows:
\begin{equation}
\label{eq:Fourier}
\mathcal{F}( (-\Delta)^s w) (\xi) = | \xi |^{2s} \mathcal{F}(w) (\xi).
\end{equation}
Equivalently, $(-\Delta)^s$ can be defined by means of the following pointwise formula
\begin{equation}
 (-\Delta)^s w(x) = C(n,s) \, \mathrm{ p.v } \int_{\mathbb{R}^n} \frac{w(x) - w(y)}{|x-y|^{n+2s}} \mathrm{d}y,
 \qquad
 C(n,s) = \frac{2^{2s} s \Gamma(s+\frac{n}{2})}{\pi^{n/2}\Gamma(1-s)},
\end{equation}
where $\textrm{p.v}$ stands for the Cauchy principal value and $C(n,s)$ is a positive normalization constant that depends only on $n$ and $s$ \cite[equation (3.2)]{NPV:12} and is introduced to guarantee that the symbol of the resulting operator is $| \xi |^{2s}$. We refer the reader to \cite[Section 1.1]{Landkof} and \cite[Proposition 3.3]{NPV:12} for a proof of the equivalence of these two definitions. Note that, as previously mentioned, there exist other non-equivalent definitions of the fractional Laplacian on bounded domains, e.g. the regional fractional Laplacian, the spectral fractional Laplacian, etc. We refer the reader to \cite{MR3393253,MR3246044,MR2270163} for a comprehensive description and study.

The rest of the paper is organized as follows.
In section \ref{sec:notation} we introduce some notation that will be useful throughout the paper.
In section \ref{sec:optimal_control} we formulate the optimal control problem for the integral fractional Laplacian with Dirichlet volume constraints.
We also prove the well-posedness of the formulation, derive optimality conditions and derive regularity estimates for the optimal variables.
Section \ref{sec:approximation} is devoted to the study of discretization techniques to solve the fractional optimal control problem.
In section \ref{subsec:fem_state} we review the a priori error analysis developed in \cite{MR3620141} for the state equation.
In section \ref{subsec:va} we propose a semidiscrete scheme for the control problem and derive a priori error estimates for the approximation of the control variable.
In section \ref{subsec:fd} we propose a fully discrete scheme for the fractional optimal control problem and derive error estimates for the approximation of the state and control variables.
In section \ref{sec:num_experiments} we report results of two-dimensional numerical tests that illustrate the theory and demonstrate the efficient solution of the discretized fractional control problem.

\section{Notation and preliminaries}\label{sec:notation}
Throughout this work $\Omega$ is an open bounded domain of $\mathbb{R}^n$ $(n \geq 1)$ with Lipschitz boundary $\partial \Omega$ that satisfies the exterior ball condition. We will denote by $\Omega^c$ the complement of $\Omega$. The relation ${\sf a} \lesssim {\sf b}$ indicates that ${\sf a} \leq C {\sf b}$ with a constant $C$ that does not depend on neither ${\sf a}$ and ${\sf b}$ but it might depend on $s$ and $\Omega$. The value of $C$ might change at each occurrence. If $\mathcal{X}$ and $\mathcal{Y}$ are normed spaces, we write $\mathcal{X}  \hookrightarrow \mathcal{Y}$ to denote that $\mathcal{X}$ is
continuously embedded in $\mathcal{Y}$.

\subsection{Function spaces}
For any $s \geq 0$, we define $H^s(\mathbb{R}^n)$, the Sobolev space of order $s$ over $\mathbb{R}^n$, by \cite[Definition 15.7]{Tartar}
\[
 H^s(\mathbb{R}^n) := \left \{ v \in L^2(\mathbb{R}^n): (1+|\xi|^2)^{s/2} \mathcal{F}(v) \in L^2(\mathbb{R}^n)\right \}.
\]
With the space $H^s(\mathbb{R}^n)$ at hand, we define $\tilde H^s(\Omega)$ as the closure of $C_0^{\infty}(\Omega)$ in $H^s(\mathbb{R}^n)$ and note that it can be equivalently characterized by \cite[Theorem 3.29]{McLean}
\begin{equation}
\tilde H^s(\Omega) = \{v|_{\Omega}: v \in H^s(\mathbb{R}^n), \textrm{ supp } v \subset \overline\Omega\}.
\end{equation}
When $\partial \Omega$ is Lipschitz $\tilde H^s(\Omega)$ is equivalent to $\mathbb{H}^s(\Omega)=[L^2(\Omega),H_0^1(\Omega)]_s$, the real interpolation between $L^2(\Omega)$ and $H_0^1(\Omega)$ for $s \in (0,1)$ and to $H^s(\Omega) \cap H_0^1(\Omega)$ for $s \in (1,3/2)$ \cite[Theorem 3.33]{McLean}. We denote by $H^{-s}(\Omega)$ the dual space of $\tilde H^s(\Omega)$ and by $\langle \cdot, \cdot \rangle$ the duality pair between these two spaces. We also define the bilinear form
\begin{equation}
\label{eq:bilinear_form}
 \mathcal{A}(v,w) = \frac{C(n,s)}{2} \iint_{\mathbb{R}^n \times \mathbb{R}^n} \frac{( v(x) - v(y) ) (w(x)-w(y))}{|x-y|^{n+2s}} \mathrm{d}x \mathrm{d}y,
\end{equation}
and denote by $\| \cdot \|_s$ the norm that $ \mathcal{A}(\cdot,\cdot)$ induces, which is just a multiple of the $H^s(\mathbb{R}^n)$-seminorm:
\[
 \| v \|_s =  \mathcal{A}(v,v)^{\frac{1}{2}} = \sqrt{\frac{C(n,s)}{2}} |v|_{H^s(\mathbb{R}^n)}.
\]

\subsection{The state equation}
Let $f \in H^{-s}(\Omega)$. The weak formulation of the state equation \eqref{eq:state_equation} reads as follows: Find $u \in \tilde H^s(\Omega)$ such that
\begin{equation}
 \label{eq:weak_state_equation}
  \mathcal{A}(u,v) = \langle f + z,v \rangle \quad \forall v \in \tilde H^{s}(\Omega).
\end{equation}
Since $\mathcal{A}$ is continuous and coercive in $\tilde H^{s}(\Omega)$, the Lax-Milgram lemma implies that problem \eqref{eq:weak_state_equation} admits a unique solution that satisfies the stability estimate
\begin{equation}
 \| u \|_s \lesssim \| f + z \|_{H^{-s}(\Omega)}.
 \label{eq:stability}
\end{equation}

\section{The fractional optimal control problem}\label{sec:optimal_control}
In this section, we analyze the \emph{fractional optimal control problem} \eqref{eq:min}--\eqref{eq:control_constraints}. We derive existence and uniqueness results together with first order necessary and sufficient optimality conditions and regularity estimates.

For $J$ defined in \eqref{eq:cost_functional}, the fractional optimal control problem reads as follows: Find $\min J(u,z)$ subject to the state equation \eqref{eq:weak_state_equation} and the control constraints \eqref{eq:control_constraints}.  The set of admissible controls is defined by
\begin{equation}
 \mathbb{Z}_{\textrm{ad}}:=
 \{  w \in L^2(\Omega): a \leq w(x) \leq b \textrm{ a.e. } x \in \Omega \},
\end{equation}
which is a nonempty, bounded, closed, and convex subset of $L^2(\Omega)$.

As it is customary in optimal control theory \cite{MR0271512,MR2583281}, to analyze \eqref{eq:min}--\eqref{eq:control_constraints}, we introduce the so-called control-to-state operator.
\begin{definition}[control-to-state map]
$\mathbf{S}: L^2(\Omega) \ni z \mapsto u(z) \in \tilde H^s(\Omega)$, where $u(z)$ solves \eqref{eq:weak_state_equation}, is called the fractional control to state operator.
\end{definition}

We notice that $\mathbf{S}$ is affine. In fact,
\begin{equation}
\mathbf{S} z = \mathbf{S}_0  z + \psi_0,
\label{eq:affine}
\end{equation}
where $\mathbf{S}_0 z$ denotes the solution to \eqref{eq:weak_state_equation} with $f \equiv 0$ and $\psi_0$ solves \eqref{eq:weak_state_equation} with $z \equiv 0$;  the operator $\mathbf{S}_0$ is linear. We also notice that $\mathbf{S}$ is self-adjoint and, in light of the estimate \eqref{eq:stability}, it is a continuous operator. In view of the continuous embeddings $H^{-s}(\Omega) \hookrightarrow L^2(\Omega) \hookrightarrow \tilde H^s(\Omega)$ \cite[Theorem 3.27]{McLean}, we may also consider $\mathbf{S}$ acting from $L^2(\Omega)$ onto itself. For simplicity, we keep the notation $\mathbf{S}$.

An optimal fractional state-control pair is defined as follows.

\begin{definition}[optimal fractional state-control pair]
A state-control pair $(\bar u(\bar z),\bar z) \in \tilde H^s(\Omega) \times \mathbb{Z}_{\mathrm{ad}}$ is called optimal for problem \eqref{eq:min}--\eqref{eq:control_constraints} if $\bar u (\bar z) = \mathbf{S} \bar z$ and
 \[
  J(\bar u (\bar z),\bar z) \leq  J(u ( z), z)
 \]
for all $(u ( z), z) \in \tilde H^s(\Omega) \times \mathbb{Z}_{\mathrm{ad}}$ such that $u(z) = \mathbf{S} z$.
\end{definition}

\smallskip The existence and uniqueness of an optimal state--control pair is as follows.

\smallskip
\begin{theorem}[existence and uniqueness]
The fractional optimal control problem \eqref{eq:min}--\eqref{eq:control_constraints} has a unique solution $(\bar u, \bar z) \in \tilde H^s(\Omega) \times \mathbb{Z}_{\textrm{ad}}$.
\end{theorem}
{\it Proof.}
By definition of $\mathbf{S}$, problem  \eqref{eq:min}--\eqref{eq:control_constraints} reduces to the following quadratic optimization problem: Minimize
\begin{equation}
\label{eq:reduced_cost_functional}
j(z) := \frac{1}{2} \| \mathbf{S} z  - u_{d} \|_{L^2(\Omega)}^2 + \frac{\alpha}{2} \| z \|_{L^2(\Omega)}^2
\end{equation}
over the set $\mathbb{Z}_{\textrm{ad}}$. Since $\alpha > 0$, it is immediate that the functional $j$ is strictly convex. In addition, since $\mathbf{S}$ is continuous, $j$ is weakly lower semicontinuous. On the other hand, the set $\mathbb{Z}_{\textrm{ad}}$ is weakly sequentially compact. The assertion thus follows from employing the direct method of the calculus of variations \cite[Theorem 5.51]{MR3026831}.
$\boxempty$

\subsection{First order optimality conditions}
To provide first order necessary and sufficient optimality conditions, we introduce the so-called adjoint state.
\begin{definition}[fractional adjoint state]
The solution $p = p(z) \in \tilde H^s(\Omega)$ of
\begin{equation}
 \label{eq:weak_adjoint_equation}
 \mathcal{A}(v,p) = \langle u - u_d , v\rangle \quad \forall v \in \tilde H^s(\Omega)
\end{equation}
is called the fractional adjoint state associated to $u = u(z)$.
\label{def:adjoint_state}
\end{definition}

The following theorem proves necessary and sufficient optimality conditions for the optimal control problem \eqref{eq:min}--\eqref{eq:control_constraints}.
\begin{theorem}[first order optimality conditions]
 $\bar z \in \mathbb{Z}_{\textrm{ad}}$ is the optimal control of problem \eqref{eq:min}--\eqref{eq:control_constraints} if and only if it satisfies the variational inequality
 \begin{equation}
  \label{eq:variational_inequality}
  (\bar p + \alpha \bar z, z - \bar z)_{L^2(\Omega)} \geq 0
 \end{equation}
 for every $z \in \mathbb{Z}_{\textrm{ad}}$, where $\bar p = \bar p (\bar z)$ solves \eqref{eq:weak_adjoint_equation} with $u$ replaced by $\bar u={\bf S}\bar z$.
\end{theorem}
{\it Proof.}
A classical result \cite[Lemma 2.21]{MR2583281} guarantees that
$\bar z \in \mathbb{Z}_{\textrm{ad}}$ minimizes the reduced cost functional $j$, defined as in \eqref{eq:reduced_cost_functional}, if and only if
\begin{equation}
(j'(\bar z) , z - \bar z)_{L^2(\Omega)} \geq 0
\label{eq:first_optimal}
\end{equation}
for every $z \in \mathbb{Z}_{\textrm{ad}}$. By standard arguments, we conclude that $j$ is Fr\'echet differentiable and we rewrite \eqref{eq:first_optimal} as
\[
(\mathbf{S} \bar z - u_d ,\mathbf{S}_0( z - \bar z))_{L^2(\Omega)} + \alpha (\bar z , z - \bar z)_{L^2(\Omega)}  \geq 0 \quad \forall z \in \mathbb{Z}_{\mathrm{ad}},
\]
where $\mathbf{S}_0$ is defined in \eqref{eq:affine} \cite[Theorem 2.20]{MR2583281}. Notice that $\mathbf{S}_0$ is self-adjoint. We can thus utilize Definition \ref{def:adjoint_state} to conclude that $\mathbf{S}_0 (\mathbf{S} \bar z - u_d) + \alpha \bar z  = \bar p + \alpha \bar z$. This concludes the proof. 
$\boxempty$

\subsection{Regularity of the optimal control}
In order to derive a priori error estimates for the solution techniques that we will propose in section \ref{subsec:va} and \ref{subsec:fd}, it is fundamental to study the regularity properties of the optimal variables associated to \eqref{eq:min}--\eqref{eq:control_constraints}. To accomplish this task, we introduce the projection operator $\textrm{proj}: L^1(\Omega) \rightarrow \mathbb{Z}_{\mathrm{ad}}$, which is defined by
\begin{equation}
 \label{eq:projection_formula}
 \textrm{proj}_{[a ,b]}(v) (x) =  \min \{ b ,\max \{ a  , v(x) \}  \} \textrm{ for all } x \in \overline \Omega,
\end{equation}
where $a$ and $b$ are in $\mathbb{R}$. With this nonlinear operator at hand, the arguments developed in \cite[Section 2.8]{MR2583281} allow us to conclude the following result: If $\alpha > 0$ and $\bar p $ is given by Definition \ref{def:adjoint_state}, then the variational inequality \eqref{eq:variational_inequality} is equivalent to the following projection formula:
\begin{equation}
 \label{eq:projection_formula_control}
 \bar z(x) = \textrm{proj}_{[a, b ]}\left( -\frac{1}{\alpha} \bar p(x) \right).
\end{equation}

\subsubsection{Regularity results on smooth domains}
We now state a regularity result for the state equation \eqref{eq:weak_state_equation} that is instrumental to derive regularity estimates for the optimal control variables.
\begin{proposition}[regularity of $u$ on smooth domains]
 Let $s \in (0,1)$ and $\Omega$ be a domain such that $\partial \Omega \in C^{\infty}$. If $f+z \in H^{r}(\Omega)$, for some $r \geq -s$, then the solution $u$ of problem \eqref{eq:weak_state_equation} belongs to $H^{s + \vartheta}(\Omega)$, where $\vartheta = \min \{s+r,1/2-\epsilon \}$ and $\epsilon > 0$ is arbitrarily small. In addition, the following estimate holds:
\begin{equation}
\| u \|_{H^{s+\vartheta}(\Omega)} \lesssim \| f + z\|_{H^{r}(\Omega)},
\label{eq:regularity_state_smooth}
 \end{equation}
where the hidden constant depends on the domain $\Omega$, $n$, $s$, and $\vartheta$.
\label{pro:state_regularity_smooth}
\end{proposition}
{\it Proof.}
 See \cite{MR3276603}.
$\boxempty$

The following example shows that, even when $\partial \Omega$ is smooth, smoothness of the right hand side $f+z$ does not ensure that solutions are any smoother than $H^{s+1/2-\epsilon}(\Omega)$ \cite{Getoor,MR3447732}: Consider $\Omega = B(0,1) \subset \mathbb{R}^n$ and $f + z \equiv 1$, then the solution to \eqref{eq:weak_state_equation} is given by
\begin{equation}
\label{eq:example}
u(x) = \frac{\Gamma(\frac{n}{2})}{2^{2s}\Gamma(\frac{n+2s}{2})\Gamma(1+s)} \left( 1 - |x|^{2}\right)^s_{+},
\end{equation}
where $t^{+} = \max \{ t,0\}$.

With the regularity estimates of Theorem \ref{thm:reg_state_smooth} at hand, we now proceed to investigate the regularity properties of the optimal control variable $\bar z$ when $\partial \Omega \in C^{\infty}$.
\begin{theorem}[regularity of $\bar z$ on smooth domains]
 Let $u_d \in H^{\lambda}(\Omega)$ with \(\lambda=\min\{1-2s,\frac{1}{2}-s-\epsilon\}\) and $f \in H^{\beta}(\Omega)$ with \(\beta=\max\{-s,\frac{1}{2}-3s-\epsilon\}\) where $\epsilon >0$ is arbitrarily small.
Then $\bar z \in H^{\gamma}(\Omega)$ with \(\gamma=\min\{1,\frac{1}{2}+s-\epsilon\}\).
In addition, we have that
\begin{equation}
\| \bar z \|_{H^{\gamma}(\Omega)} \lesssim \| f \|_{H^{\beta}(\Omega)} + 
\| \bar z\|_{L^{2}(\Omega)} +  \| u_d \|_{H^{\lambda}(\Omega)},
\label{eq:reg_control_smooth2}
\end{equation}
where the hidden constant depends on $\Omega$, $n$, and $s$.
\label{thm:reg_control_smooth}
\end{theorem}
{\it Proof.}
We begin by noticing that, since the right-hand sides of the state and adjoint equations, namely, $f + \bar z$ and $\bar u - u _d$, respectively, belong to $H^{-s}(\Omega)$, we have that $\bar u, \bar p \in H^s(\Omega)$.
This, on the basis of a nonlinear operator interpolation result as in \cite[Lemma 28.1]{Tartar} combined with \cite[Theorem A.1]{MR1786735} and formula \eqref{eq:projection_formula_control}, implies that $\bar z \in H^{s}(\Omega)$.

We now consider the following cases.

\smallskip

\noindent \emph{Case} 1. \(s\geq 1/4\): Notice that, in view of the assumption on $u_d$, we have that \(\bar u - u_{d}\in H^{\eta}(\Omega)\), where \(\eta=\min\{s,\lambda\}=\lambda\). By Proposition~\ref{pro:state_regularity_smooth}, we conclude that \(\bar p\in H^{\xi}(\Omega)\), where \(\xi= s + \vartheta_1\) and $\vartheta_1 = \min\{s+\lambda,1/2-\epsilon\}$. By invoking, again, \cite[Lemma 28.1]{Tartar}, \cite[Theorem A.1]{MR1786735}, and formula \eqref{eq:projection_formula_control}, we conclude that \(\bar z\in H^{\kappa}(\Omega)\) with \(\kappa=\min\{1,\xi\}\).

Notice that, if $s \in [1/4,1/2 + \epsilon)$, we have that $1 - 2s >1/2 - s - \epsilon$ and thus that $\lambda = 1/2 - s - \epsilon$. Consequently, $\vartheta_1 = 1/2 - \epsilon$ and $\xi = s + \vartheta_1 = s + 1/2 - \epsilon$. As a result, we have obtained that $\bar z \in H^{s + 1/2 - \epsilon}(\Omega)$ with $\epsilon > 0$ being arbitrarily small. On the other hand, if $s \geq 1/2 + \epsilon$, then $1 - 2s \leq 1/2 - s - \epsilon$. This yields $\lambda = 1-2s$. Consequently, $\vartheta_1 = \min \{ 1-s, 1/2 - \epsilon \} = 1-s$, which implies that $\xi = s + \vartheta_1 = 1$. We have thus obtained that $\bar z \in H^1(\Omega)$.

\smallskip

\noindent \emph{Case} 2. \(1/8 \leq s < 1/4\): Notice that, since \(\beta=1/2-3s-\epsilon<s\), we obtain that \(f + \bar z \in H^{\beta}(\Omega)\). We can thus apply Proposition~\ref{pro:state_regularity_smooth} to conclude that  \(\bar u\in H^{\iota}(\Omega)\) with $\iota = s + \vartheta_2$ and $\vartheta_2 = \min \{s + \beta ,1/2 - \epsilon\}$. Notice that $\iota = 1/2-s-\epsilon = \lambda$.
Since \(u_{d}\in H^{\lambda}(\Omega)\), we thus have that \(\bar u - u_{d}\in H^{\lambda}(\Omega)\). Therefore, by Proposition~\ref{pro:state_regularity_smooth}, $\bar p \in H^{s + \vartheta_3}(\Omega)$, where $\vartheta_3 = \min \{ s + \lambda, 1/2 - \epsilon \} = 1/2 - \epsilon$. By invoking, again, \cite[Lemma 28.1]{Tartar}, \cite[Theorem A.1]{MR1786735}, and formula \eqref{eq:projection_formula_control}, we obtain that $\bar z \in H^{s + 1/2 - \epsilon}(\Omega)$, where $\epsilon >0$ is arbitrarily small.

\smallskip

\noindent \emph{Case} 3. \(0 < s < 1/8\): Since $\bar z \in H^s(\Omega)$, we have that $f + \bar z \in H^{\delta}(\Omega)$, where $\delta = \min \{ s, 1/2-3s-\epsilon\} = s$. We thus invoke Proposition~\ref{pro:state_regularity_smooth} to conclude that \(\bar u \in H^{s + \vartheta_4}(\Omega)\), where \(\vartheta_4 = \min \{ s + \delta, 1/2 - \epsilon\} = \min \{ 2s, 1/2 - \epsilon \} = 2s\).  Notice that, in view of the assumption \(u_{d}\in H^{\lambda}(\Omega)\) with \(\lambda=1/2-s-\epsilon\), we conclude that \(\bar u-u_{d} \in H^{2s}(\Omega)\). We apply again \cite[Lemma 28.1]{Tartar}, \cite[Theorem A.1]{MR1786735}, and formula \eqref{eq:projection_formula_control} to conclude that  \(\bar z\in H^{s + \vartheta_5}(\Omega)\), where $\vartheta_5 = \min \{ 3s, 1/2 - \epsilon \}$.

\smallskip

\hspace{0.1cm}
\noindent \emph{Case} 3.1.  \(1/6 \leq s < 1/8\): In this case $3s > 1/2 - \epsilon$, and then $\vartheta_5 = 1/2 - \epsilon$. Consequently, $\bar z \in H^{s + 1/2 - \epsilon}(\Omega)$ with $\epsilon > 0$ being arbitrarily small.

\hspace{0.1cm}
\noindent \emph{Case} 3.2.  \(0 < s < 1/6\): On the basis of the arguments previously developed, a bootstrap argument allows us to conclude that $\bar z \in H^{s + 1/2 - \epsilon}(\Omega)$ with $\epsilon > 0$ being arbitrarily small.

In all the considered cases, the estimate \ref{eq:reg_control_smooth2} follows from stability estimates for state and adjoint equations and the nonlinear operator interpolation result of \cite[Lemma 28.1]{Tartar} combined with \cite[Theorem A.1]{MR1786735} and formula \eqref{eq:projection_formula_control}. This concludes the proof.
$\boxempty$

The following result follows immediately.
\begin{corollary}[regularity of $\bar u$ and $\bar p$ on smooth domains]
Let $s \in (0,1)$. Under the framework of Theorem \ref{thm:reg_control_smooth} we have that $\bar u \in H^{s + 1/2 - \epsilon}(\Omega)$ and $\bar p \in H^{s + 1/2 - \epsilon}(\Omega)$ for every $\epsilon >0$.
\label{thm:reg_state_smooth}
\end{corollary}

As the example previously described, which involves \eqref{eq:example} as exact solution, shows, the regularity properties of the optimal variables $\bar u$ and $\bar p$ obtained in Corollary \ref{thm:reg_state_smooth} cannot be improved.

\subsubsection{Regularity results on Lipschitz domains}
The following results establish regularity estimates in H{\"o}lder and Sobolev spaces for Lipschitz domains.
\begin{proposition}[regularity of $u$ on Lipschitz domains]
Let $s \in (0,1)$ and $\Omega$ be a bounded Lipschitz domain satisfying the exterior ball condition. If $f + z  \in L^{\infty}(\Omega)$, then the solution $u$ of problem \eqref{eq:weak_state_equation} belongs to $C^s(\mathbb{R}^n)$ and the following estimate holds:
 \begin{equation}
\| u \|_{C^s(\mathbb{R}^n)} \lesssim \| f \|_{L^{\infty}(\Omega)} + \| z \|_{L^{\infty}(\Omega)},
\label{eq:regularity_state}
 \end{equation}
where the hidden constant depends on $\Omega$ and $s$.
\label{pro:state_regularity_1}
\end{proposition}
{\it Proof.}
 See \cite[Proposition 1.1]{MR3168912}.
$\boxempty$
\begin{proposition}[regularity of $u$ on Lipschitz domains]
Let $s \in (0,1)$ and $\Omega$  be a bounded Lipschitz domain satisfying the exterior ball condition. If $s \in (0,1/2)$, let $f+z \in C^{\frac{1}{2}-s}(\overline \Omega)$; if $s = 1/2$, let $f+z \in L^{\infty}(\Omega)$; and if $s \in (1/2,1)$, let $f +z \in C^{\beta}(\overline \Omega)$ for some $\beta >0$. Then, for every $\epsilon >0$, the solution $u$ of problem \eqref{eq:weak_state_equation} belongs to $H^{s+1/2-\epsilon}(\Omega)$ and satisfies the estimate
 \begin{equation}
\| u \|_{H^{s+1/2-\epsilon}(\Omega)} \lesssim \| f + z \|_{\star},
\label{eq:regularity_state_12}
 \end{equation}
where $\|  \cdot \|_{\star}$ denotes the $C^{\frac{1}{2}-s}(\overline \Omega)$, $L^{\infty}(\Omega)$ or $C^{\beta}(\Omega)$-norm, correspondingly to whether $s$ is smaller, equal or grater than $1/2$. The hidden constant depends on the domain $\Omega$, the dimension $n$, and the parameter $s$, and blows up when $\epsilon \rightarrow 0$.
\label{pro:state_regularity_12}
\end{proposition}
{\it Proof.}
 See \cite[Propositions 3.6 and 3.11]{MR3620141}.
$\boxempty$

We now proceed to investigate the regularity properties of the optimal control variable $\bar z$ when $\Omega$ is a  bounded Lipschitz domain that satisfies the exterior ball condition. We begin with the case $s \in (0,\tfrac{1}{4})$.
\begin{theorem}[regularity of $\bar z$ on Lipschitz domains: $s \in (0,\tfrac{1}{4})$]
 Let $f \in L^{\infty}(\Omega)$ and $u_d \in L^{\infty}(\Omega)$. If $s \in (0,\tfrac{1}{4})$, then we have that $\bar z \in C^{s}(\overline \Omega)$. In addition, we have the estimate
\[
\| \bar z \|_{C^{s}(\overline \Omega)} \lesssim \| f \|_{L^{\infty}(\Omega)} + \| \bar z \|_{L^{\infty}(\Omega)} + \| u_d \|_{L^{\infty}(\Omega)},
\]
where the hidden constant depends on $\Omega$ and $s$.
\label{thm:reg_control_0}
\end{theorem}
{\it Proof.}
 Since the right-hand side $f + \bar z$ of the state equation \eqref{eq:weak_state_equation} belongs to $L^{\infty}(\Omega)$, Proposition \ref{pro:state_regularity_1} allows us to conclude that $\bar u \in C^s(\mathbb{R}^n)$. Thus, since $u_d \in L^{\infty}(\Omega)$, we can apply Proposition \ref{pro:state_regularity_1}, again, to conclude that $\bar p \in C^s(\overline \Omega)$. The projection formula \eqref{eq:projection_formula} and \cite[Theorem A.1]{MR1786735} allow us to conclude that $\bar z \in C^s(\overline \Omega)$.
$\boxempty$
\begin{theorem}[regularity of $\bar z$ on Lipschitz domains: $s \in [\frac{1}{4},\frac{1}{2})$]
Let $f \in L^{\infty}(\Omega)$ and $u_d \in C^{1/2-s}(\overline \Omega)$. If $s \in [\tfrac{1}{4},\tfrac{1}{2})$, then we have that, for every $\epsilon >0$, the optimal control $\bar z \in H^{s + 1/2 -\epsilon}(\Omega)$. In addition, we have the estimate
\[
\| \bar z \|_{H^{s + 1/2 -\epsilon}(\Omega)} \lesssim\| f \|_{L^{\infty}(\Omega)} + \| \bar z \|_{L^{\infty}(\Omega)} + \| u_d \|_{C^{1/2-s}(\overline \Omega)},
\]
where the hidden constant depends on $\Omega$, $n$, and $s$, and blows up when $\epsilon \rightarrow 0$.
\label{thm:reg_control_1}
\end{theorem}
{\it Proof.}
In view of the fact that $f + \bar z$ belongs to $L^{\infty}(\Omega)$, we can apply the results of Proposition \ref{pro:state_regularity_1} to obtain that $\bar u \in C^s(\mathbb{R}^n)$ and that
\begin{equation}
 \| \bar u \|_{C^s(\mathbb{R}^n)} \lesssim \| f \|_{L^{\infty}(\Omega)} + \| \bar z \|_{L^{\infty}(\Omega)}.
 \label{eq:u_Cs_1}
\end{equation}
Now, notice that, since $s \in [1/4,1/2)$, the following trivial inequality holds: $1/2-s \leq s$. This, the estimate  \eqref{eq:u_Cs_1}, and the assumption on the desired state $u_d$ reveal that $\bar u - u_d \in C^{1/2-s}(\overline \Omega)$. We are thus in position to apply the results of Proposition \ref{pro:state_regularity_12} to obtain that, for every $\epsilon >0$, the optimal adjoint variable $\bar p$ belongs to $H^{s+1/2-\epsilon}(\Omega)$. In addition, we have the estimate
\[
 \| \bar p \|_{H^{s+1/2-\epsilon}(\Omega)} \lesssim \| \bar u \|_{C^{1/2-s}(\overline \Omega)} + \| u_d \|_{C^{1/2-s}(\overline \Omega)}.
\]
In view of the projection formula \eqref{eq:projection_formula} and \cite[Theorem A.1]{MR1786735}, a nonlinear operator interpolation result as in \cite[Lemma 28.1]{Tartar} allow us to conclude that, for every $\epsilon > 0$, $\bar z \in H^{s + 1/2 -\epsilon}(\Omega)$, with the estimate
\[
\| \bar z \|_{H^{s + 1/2 -\epsilon}(\Omega)} \lesssim \| \bar p \|_{H^{s+1/2-\epsilon}(\Omega)} \lesssim \| \bar u \|_{C^{1/2-s}(\overline \Omega)} + \| u_d \|_{C^{1/2-s}(\overline \Omega)}.
\]
This, in view of \eqref{eq:u_Cs_1}, concludes the proof.
$\boxempty$

We now consider the case $s \in ( \frac{1}{2},1)$.

\begin{theorem}[regularity of $\bar z$ on Lipschitz domains: $s \in (\tfrac{1}{2},1)$]
Let $f \in L^{\infty}(\Omega)$ and $u_d \in C^{\beta}(\Omega)$, for some $\beta > 0$. If $s \in (\tfrac{1}{2},1)$, then we have that the optimal control $\bar z$ belongs to $H^{1}(\Omega)$. In addition, we have the estimate
\[
 \| \bar z \|_{H^{1}(\Omega)} \lesssim
 \| f \|_{L^{\infty}(\Omega)} + \| \bar z \|_{L^{\infty}(\Omega)} + \| u_d \|_{C^{\gamma}(\overline \Omega)},
\]
where $\gamma = \min \{\beta,s\}$, and the hidden constant depends on $\Omega$, $n$, and $s$, and blows up when $\epsilon \rightarrow 0$.
\label{thm:reg_control_2}
\end{theorem}
{\it Proof.}
We begin the proof by applying the results of Proposition \ref{pro:state_regularity_1} to conclude that $\bar u \in C^s(\mathbb{R}^n)$, with the estimate
\begin{equation}
 \| \bar u \|_{C^s(\mathbb{R}^n)} \lesssim \| f \|_{L^{\infty}(\Omega)} + \| \bar z \|_{L^{\infty}(\Omega)}.
 \label{eq:aux_u_Cs}
\end{equation}
In view of the assumptions, we conclude that $\bar u - u_d \in C^{\gamma}(\overline \Omega)$, where $\gamma = \min \{\beta,s\}$. We can thus invoke the results of Proposition \ref{pro:state_regularity_12} to conclude that, for every $\epsilon >0$, we have that $\bar p \in H^{s+1/2-\epsilon}(\Omega)$, with the estimate
\[
 \| \bar p \|_{H^{s+1/2-\epsilon}(\Omega)} \lesssim \| \bar u \|_{C^{\gamma}(\overline \Omega)} + \| u_d \|_{C^{\gamma}(\overline \Omega)}.
\]
The regularity property for the optimal control follows thus from \eqref{eq:projection_formula}, \cite[Theorem A.1]{MR1786735} and \cite[Lemma 28.1]{Tartar}. In fact, we have that $\bar z \in H^{1}(\Omega)$, with the estimate
\begin{align*}
\| \bar z \|_{H^{1}(\Omega)} & \lesssim \| \bar u \|_{C^{\gamma}(\overline \Omega)} + \| u_d \|_{C^{\gamma}(\overline \Omega)}
\\
& \lesssim  \| f \|_{L^{\infty}(\Omega)} + \| \bar z \|_{L^{\infty}(\Omega)} + \| u_d \|_{C^{\gamma}(\overline \Omega)},
\end{align*}
where, to obtain the last estimate, we have used \eqref{eq:aux_u_Cs}. This concludes the proof.
$\boxempty$

Similar arguments to the ones elaborated in the proofs of Theorems \ref{thm:reg_control_1} and \ref{thm:reg_control_2} allow us to obtain regularity estimates for the case $s=\frac{1}{2}$. For brevity, we present the following result and skip the details.
\begin{theorem}[regularity of $\bar z$ on Lipschitz domains: $s = \tfrac{1}{2}$]
Let $f$ and $u_d \in L^{\infty}(\Omega)$. If $s=  \tfrac{1}{2}$, then we have that, for every $\epsilon >0$, the optimal control $\bar z \in H^{1-\epsilon}(\Omega)$, with the estimate
\[
 \| \bar z \|_{H^{1-\epsilon}(\Omega)} \lesssim  \| f \|_{L^{\infty}(\Omega)} + \| \bar z \|_{L^{\infty}(\Omega)} + \| u_d \|_{L^{\infty}(\Omega)},
\]
where the hidden constant depends on $\Omega$, $n$, and $s$, and blows up when $\epsilon \rightarrow 0$.
\label{thm:reg_control_3}
\end{theorem}

The following regularity result will be instrumental for the error analysis that we will perform.
\begin{lemma}[regularity of $\bar z$ on Lipschitz domains: $s \in [\tfrac{1}{4},1)$]
Let $f \in L^{\infty}(\Omega)$ and $u_d \in L^{\infty}(\Omega)$. In addition, for $s \in [1/4,1/2)$, let $u_d \in C^{\beta}(\Omega)$ for some $\beta > 0$. Then,
\begin{equation}
\label{eq:regularity_z_as_rhs}
\bar z \in
\begin{dcases}
C^{1/2-s}(\overline \Omega), & s \in [\tfrac{1}{4},\tfrac{1}{2}), \\
L^{\infty}(\Omega), & s = \tfrac{1}{2}, \\
C^{s}(\overline \Omega), & s \in (\tfrac{1}{2},1).
\end{dcases}
\end{equation}
\end{lemma}
{\it Proof.}
The case $s = \tfrac{1}{2}$ follows immediately from the fact that $\bar z \in \mathbb{Z}_{\mathrm{ad}}$.

If $s \in (\tfrac{1}{2},1)$, we can apply Proposition \ref{pro:state_regularity_1}, since $\bar u - u_d \in L^{\infty}(\Omega)$, to conclude that $\bar p \in C^s(\mathbb{R}^n)$. This, in view of the projection formula \eqref{eq:projection_formula} reveals that $z \in C^s(\overline \Omega)$.

If $s \in [\tfrac{1}{4},\tfrac{1}{2})$, an application of Proposition \ref{pro:state_regularity_1}, again, yields $\bar p \in C^{s}(\overline \Omega)$. This implies that $\bar p \in C^{1/2-s}(\overline \Omega)$ for $s \in [\tfrac{1}{4},\tfrac{1}{2})$. The projection formula \eqref{eq:projection_formula} allows us to conclude.
$\boxempty$

\section{Approximation of the fractional control problem}\label{sec:approximation}
In this section, we introduce and analyze two solution techniques to approximate the solution to the fractional optimal control problem \eqref{eq:min}--\eqref{eq:control_constraints}. Before proceeding with the design and analysis of the proposed methods, it is instructive to review the numerical approximation of the state equation \eqref{eq:state_equation} developed in \cite{MR3620141}. We briefly report such results in the following section.

\subsection{A finite element method for the state equation}
\label{subsec:fem_state}
We start with some terminology and describe the construction of the underlying finite element spaces. Let $\T = \{ T \}$ be a conforming partition of $\overline \Omega$ into simplices $T$ with size $h_T = \diam(T)$, and set $h_{\T} = \max_{T \in \T} h_T$. We denote by $\mathbb{T}$ the collection of conforming and shape regular meshes that are refinements of an initial mesh $\T_0$. By shape regular we mean that there exists a constant $\sigma > 1$ such that $\max \{ \sigma_T: T \in \T \} \leq \sigma$ for all $\T \in \mathbb{T}$. Here $\sigma_T = h_T / \rho_T$ denotes the shape coefficient of
$T$, where $\rho_T$ is the diameter of the largest ball that can be inscribed in $T$ \cite{MR2373954,MR1930132,MR2050138}.

Given a mesh $\T \in \mathbb{T}$, we define the finite element space of continuous piecewise polynomials of degree one as
\begin{equation}
\V(\T) = \left\{ v_{\T} \in C^0( \overline\Omega): {v_{\T}}_{|T} \in \mathbb{P}_1(T) \ \forall T \in \T, \ v_{\T} = 0 \textrm{ on } \partial \Omega \right\}.
\label{eq:defFESpace}
\end{equation}
Note that discrete functions are trivially extended by zero to $\Omega^c$ and that we enforce a classical homogeneous Dirichlet boundary condition at the degrees of freedom that are located at the boundary of \(\Omega\).
  As Proposition~\ref{pro:state_regularity_1} states, the solutions of state and adjoint equations are in the H\"older space \(C^{s}\left(\mathbb{R}^{n}\right)\). Therefore their boundary trace is zero on \(\partial \Omega\).
The finite element approximation of the state equation \eqref{eq:weak_state_equation} is then the unique solution to the following discrete problem: Find $u_{\T} \in \V(\T)$ such that
\begin{equation}
\label{eq:weak_state_discrete}
\mathcal{A}(u_{\T},v_{\T}) = \langle f + z, v_{\T} \rangle \quad \forall v_{\T} \in \V(\T),
\end{equation}
Note that discrete functions are trivially extended by zero to $\Omega^c$.
From this formulation it follows that $u_{\T}$ is the projection (in the energy
norm) of $u$ onto $\V(\T)$. Consequently, we have a C\'ea-like best approximation result
\begin{equation}
 \| u - u_{\T} \|_s =  \inf_{v_\T \in \V(\T)} \| u - v_{\T} \|_s.
\end{equation}

\subsubsection{Error estimates on quasi-uniform meshes}
Localization results for fractional seminorms \cite{MR1930387} and local stability and approximation properties for the Scott-Zhang interpolation operator \cite{MR3118443} are the key ingredients to provide an a priori error analysis. We present the following a priori error estimate in energy norm \cite[Theorem 4.7]{MR3620141}.
\begin{proposition}[energy error estimate for quasi--uniform meshes]\label{prop:HsstateQuasiUniform}
Let $u \in \tilde H^s(\Omega)$ be the solution to \eqref{eq:weak_state_equation}, and let $u_{\T} \in \V(\T)$ be the solution to the discrete problem \eqref{eq:weak_state_discrete}. If $\T$ is quasi--uniform, then, under the hypotheses of Proposition \ref{pro:state_regularity_12}, we have the error estimate
\begin{equation}
 \| u - u_{\T} \|_s \lesssim h_{\T}^{\frac{1}{2}} | \log h_{\T}|  \|  f + z \|_{\star},
 \label{eq:energy_estimate}
\end{equation}
where the hidden constant depends on $\Omega$, $s$, and $\sigma$; $\|  \cdot \|_{\star}$ denotes the $C^{\frac{1}{2}-s}(\overline \Omega)$, $L^{\infty}(\Omega)$ or $C^{\beta}(\Omega)$-norm, correspondingly to whether $s$ is smaller, equal or grater than $1/2$.
\end{proposition}
The following a priori error estimate in $L^2(\Omega)$ can be derived following the arguments of \cite[Proposition 4.3]{Borthagaray2018}; see \cite[Proposition 3.8]{CVS}.
\begin{proposition}[$L^2$-error estimate for quasi--uniform meshes]\label{prop:L2stateQuasiUniform}
Let $u \in \tilde H^s(\Omega)$ be the solution to \eqref{eq:weak_state_equation}, and let $u_{\T} \in \V(\T)$ be the solution to the discrete problem \eqref{eq:weak_state_discrete}. If $\T$ is quasi--uniform, then, under the hypotheses of Proposition \ref{pro:state_regularity_smooth}, we have the error estimate
\begin{equation}
 \| u - u_{\T} \|_{L^2(\Omega)} \lesssim h_{\T}^{\vartheta + \beta}  \|  f + z \|_{H^r(\Omega)},
 \label{eq:L2estimater}
\end{equation}
where $\vartheta =\min \{ s + r,1/2 - \epsilon \}$, $\beta = \min \{ s,1/2 - \epsilon \}$ and $\epsilon >0$ may be taken arbitrarily small. In addition, the hidden constant depends on $\Omega$, $s$, $n$, $\vartheta$, and $\sigma$ and blows up when $\epsilon \rightarrow 0$.
\end{proposition}

\subsubsection{Error estimates on graded meshes}
When $s \in (1/2,1)$ and $n=2$, the singular behavior of the solution exhibited by the regularity estimates in weighted Sobolev spaces of \cite{MR3620141} can be compensated by using a priori adapted meshes. The latter, that are graded near the boundary of the domain and allow for an improvement on the priori error estimate \eqref{eq:energy_estimate}, are constructed as follows. In addition to shape regularity, we assume that the meshes $\T$ have the following property: Given a mesh parameter $h>0$ and $\mu \in [1,2]$ every element $T \in \T$ satisfies
\begin{equation}
h_T \approx C(\sigma) h^{\mu} \textrm{ if } T \cap \partial \Omega \neq \emptyset,
\quad
h_T \approx C(\sigma) h \mathrm{dist}(T,\partial \Omega)^{(\mu-1)/\mu} \textrm{ if } T \cap \partial \Omega = \emptyset,
\label{eq:graded_meshes}
\end{equation}
where $C(\sigma)$ depends only on the shape regularity constant $\sigma$ of the mesh $\T$. We notice that $\mu$ relates the mesh parameter $h$ to the number of degrees of freedom, $N$, as follows:
\begin{equation}
N \approx h_{\T}^{-2}\textrm{ if } \mu \in (1,2), \quad
N \approx h_{\T}^{-2}|\log h_{\T}|\textrm{ if } \mu = 2.
\end{equation}
The optimal choice for the parameter is $\mu =2$ and the following error estimate can be derived \cite[Theorem 4.11]{MR3620141}.
\begin{proposition}[energy error estimate for graded meshes]
Let $\Omega \subset \mathbb{R}^2$ and $s \in (1/2,1)$. Let $u \in \tilde H^s(\Omega)$ be the solution to \eqref{eq:weak_state_equation}, and let $u_{\T} \in \V(\T)$ be the solution to the discrete problem \eqref{eq:weak_state_discrete}. If $\T$ satisfies \eqref{eq:graded_meshes} with $\mu=2$ and $f + z \in C^{1-s}(\overline \Omega)$ then, we have the error estimate
\begin{equation}
 \| u - u_{\T} \|_{s} \lesssim |\log N|  N^{-\frac{1}{2}} \|  f + z \|_{C^{1-s}(\overline \Omega)},
 \label{eq:energy_estimate_graded}
\end{equation}
where the hidden constant depends on $\sigma$ and blows up when $s \rightarrow 1/2$.
\end{proposition}

\subsection{A semidiscrete scheme: the variational approach}
\label{subsec:va}
In this section, we propose a semidiscrete scheme for the fractional optimal control problem that is based on the so-called variational discretization approach. This approach, that was introduced by Hinze in \cite{Hinze:05}, discretizes only the state space; the control space $\mathbb{Z}_{\mathrm{ad}}$ is not discretized. The scheme induces a discretization of the optimal control variable by projecting the optimal discrete adjoint state into the admissible control set.

The aforementioned semidiscrete scheme reads as follows: Find $\min J(u_{\T},g)$ subject to the discrete state equation
\begin{equation}
\label{eq:weak_state_discrete_va}
\mathcal{A}(u_{\T},v_{\T}) = \langle f + g, v_{\T} \rangle \quad \forall v_{\T} \in \V(\T),
\end{equation}
and the control constraints $g \in \mathbb{Z}_{\mathrm{ad}}$. For notational convenience, we will refer to the previously defined problem as the \emph{semidiscrete optimal control problem}.

To perform an error analysis, we introduce the control-to-state operator $\mathbf{S}_{\T}: \mathbb{Z}_{\mathrm{ad}} \ni g \mapsto u_{\T} \in \V(\T)$ where $\mathbf{S}_{\T}g = u_{\T}(g)$ solves \eqref{eq:weak_state_discrete_va}. We notice that $\mathbf{S}_{\T}$ is an affine and continuous operator. In fact,
$
\mathbf{S}_{\T} g = \mathbf{S}_{\T \!,0} \,g + \psi_{\T},
$
where $\mathbf{S}_{\T\!,0} g$ denotes the solution to \eqref{eq:weak_state_discrete_va} with $f \equiv 0$ and $\psi_{\T}$ solves \eqref{eq:weak_state_discrete_va} with $g \equiv 0$;  $\mathbf{S}_{\T\!,0}$ is a linear and continuous operator.

As in section \ref{sec:optimal_control}, we denote by $(\bar u_{\T}, \bar g) \in \V(\T)  \times \mathbb{Z}_{\mathrm{ad}}$ an optimal pair solving the semidiscrete optimal control problem.

We now state the existence and uniqueness results together with first order optimality conditions.
\begin{theorem}[existence, uniqueness and optimality conditions]
The semidiscrete optimal control problem has a unique optimal solution $(\bar u_{\T},\bar g) \in \V(\T) \times  \mathbb{Z}_{\mathrm{ad}}$. In addition, the first order optimality condition
\begin{equation}
\label{eq:variational_inequality_va}
(\bar p_{\T} + \alpha \bar g, g - \bar  g )_{L^2(\Omega)} \geq 0 \quad \forall  g \in  \mathbb{Z}_{\mathrm{ad}}
\end{equation}
is necessary and sufficient.
\end{theorem}
{\it Proof.}
The proof follows standard arguments \cite{MR2583281}. For brevity, we skip the details.
$\boxempty$

We define the optimal adjoint state $\bar p_{\T} = \bar p_{\T} (\bar g)$ as the solution to
\begin{equation}
\label{eq:weak_adjoint_discrete_va}
\mathcal{A}(v_{\T},\bar p_{\T}) = \langle \bar u_{\T} - u_d, v_{\T} \rangle \quad \forall v_{\T} \in \V(\T).
\end{equation}

With these ingredients at hand, we proceed to derive an a priori error analysis for the semidiscrete optimal control problem. The proof is inspired by the arguments developed by Hinze in \cite{MR2122182}. Since, in our case, the optimal control and state variables exhibit reduced regularity properties, that are dictated by Theorem \ref{thm:reg_control_smooth} and Corollary \ref{thm:reg_state_smooth}, we present a detailed proof.
\begin{theorem}[variational approach: error estimate]
Let $s \in (0,1)$ and $u_d \in H^{1/2 - s - \epsilon}(\Omega)$, for every $\epsilon >0$. Let $(\bar u, \bar z)$ and $(\bar u_\T, \bar g)$ be the solutions to the continuous and semidiscrete optimal control problems, respectively. If $\T$ is quasi--uniform, then, under the framework of Theorem \ref{thm:reg_control_smooth}, we have the error estimate
\begin{equation}
\label{eq:error_estimate_control_va}
 \| \bar z - \bar g  \|_{L^2(\Omega)} \lesssim h_{\T}^{1/2 + \beta -\epsilon} \left( \| \bar u\|_{H^{s+1/2-\epsilon}(\Omega)} +  \| u_d \|_{H^{1/2 - s - \epsilon}(\Omega)}  + \| \bar z\|_{H^{\gamma}(\Omega)}\right),
\end{equation}
where $\beta = \min \{ s, 1/2 - \epsilon \}$, $\gamma = \min \{ s+1/2-\epsilon,1\}$, and $\epsilon > 0$ is arbitrarily small. The hidden constant depends on $\Omega$, $s$, and $n$ and blows up when $\epsilon \rightarrow 0$.
\end{theorem}
{\it Proof.}
Set $z = \bar g$ and $g = \bar z$ in the variational inequalities \eqref{eq:variational_inequality} and \eqref{eq:variational_inequality_va}, respectively and add the obtained inequalities to arrive at the estimate
\begin{equation}
 \alpha \| \bar z - \bar g  \|^2_{L^2(\Omega)} \leq (\bar p - \bar p_{\T}, \bar g - \bar z)_{L^2(\Omega)}.
 \label{eq:va_basic_estimate}
\end{equation}

We now write $\bar p = \bar p(\bar z) = \mathbf{S}_0 (\mathbf{S} \bar z - u_d )$ and $\bar p_{\T} = \bar p_{\T}(\bar g) = \mathbf{S}_{\T \! ,0} (\mathbf{S}_{\T} \bar g - u_d )$, where $\mathbf{S}$ and $\mathbf{S}_{\T}$ denote the continuous and semidiscrete control-to-state maps, respectively. With these relations at hand we can thus rewrite the estimate \eqref{eq:va_basic_estimate} as
\begin{equation*}
 \alpha \| \bar z - \bar g  \|^2_{L^2(\Omega)} \leq ( \mathbf{S}_0 (\mathbf{S} \bar z - u_d ) - \mathbf{S}_{\T\!,0} (\mathbf{S}_{\T} \bar g - u_d ), \bar g - \bar z)_{L^2(\Omega)}.
 \label{eq:va_basic_estimate2}
\end{equation*}
Adding and subtracting the term $\mathbf{S}_{\T\!,0} \mathbf{S} \bar z$, we obtain that
\begin{equation*}
 \alpha \| \bar z - \bar g  \|^2_{L^2(\Omega)} \leq ( (\mathbf{S}_0 - \mathbf{S}_{\T\!,0}) \mathbf{S} \bar z +\mathbf{S}_{\T\!,0} \mathbf{S} \bar z -\mathbf{S}_{\T,0}\mathbf{S}_{\T} \bar g +  (\mathbf{S}_{\T\!,0} - \mathbf{S}_0) u_d, \bar g - \bar z)_{L^2(\Omega)}.
 \label{eq:va_basic_estimate3}
\end{equation*}
We now add and subtract $\mathbf{S}_{\T\!,0}\mathbf{S}_{\T} \bar z$ to conclude that
\begin{multline}
 \alpha \| \bar z - \bar g  \|^2_{L^2(\Omega)} \leq ( (\mathbf{S}_0 - \mathbf{S}_{\T\!,0}) \mathbf{S} \bar z, \bar g - \bar z)_{L^2(\Omega)} +(\mathbf{S}_{\T\!,0} (\mathbf{S} - \mathbf{S}_{\T}) \bar z,\bar g - \bar z)_{L^2(\Omega)}
\\
  (\mathbf{S}_{\T\!,0}\mathbf{S}_{\T} (\bar z- \bar g),\bar g - \bar z)_{L^2(\Omega)}  +  ((\mathbf{S}_{\T\!,0} - \mathbf{S}_0) u_d, \bar g - \bar z)_{L^2(\Omega)} =: \mathrm{I} +  \mathrm{II} +  \mathrm{III} +  \mathrm{IV}.
   \label{eq:va_basic_estimate_va}
 \end{multline}

Thus, it suffices to control the terms $\mathrm{I}$, $\mathrm{II}$, $\mathrm{III}$, and $\mathrm{IV}$. We begin with the control of $\mathrm{I}$. To accomplish this task, we first notice that, since $\mathbf{S} \bar z = \bar u$, Corollary \ref{thm:reg_state_smooth} implies that $\mathbf{S} \bar z \in H^{s+1/2-\epsilon}(\Omega)$ for every $\epsilon > 0$. We can thus invoke the error estimate \eqref{eq:L2estimater} with $r = s + 1/2 - \epsilon$ to conclude that
 \[
|\mathrm{I} | \lesssim h_{\T}^{\vartheta_1 + \beta} \| \bar u\|_{H^{s+1/2-\epsilon}(\Omega)} \| \bar z - \bar g  \|_{L^2(\Omega)},
 \]
where $\vartheta_1 = \min \{2s+1/2-\epsilon, 1/2 - \epsilon\}$ and $\beta = \min \{s,1/2-\epsilon\}$. We notice that $\vartheta_1 = 1/2 - \epsilon$, and thus that $\vartheta_1 + \beta = 1/2 + \beta - \epsilon$. The control of the term $\mathrm{IV}$ follows exactly the same arguments upon exploiting the assumption $u_d \in H^{1/2 - s -\epsilon}(\Omega)$. To estimate $\mathrm{II}$, we follow similar arguments and use the continuity of the discrete operator $\mathbf{S}_{\T}$. Finally, we control the term $\mathrm{III}$ as follows:
\[
 \mathrm{III} = (\mathbf{S}_{\T} (\bar z- \bar g),\mathbf{S}_{\T,0}(\bar g - \bar z))_{L^2(\Omega)} = - \| \mathbf{S}_{\T,0}(\bar g - \bar z) \|^2_{L^2(\Omega)} \leq 0.
\]

The desired estimate \eqref{eq:error_estimate_control_va} follows from replacing the estimates we obtained for $\mathrm{I}$, $\mathrm{II}$, $\mathrm{III}$, and $\mathrm{IV}$ into \eqref{eq:va_basic_estimate_va}. This concludes the proof.
$\boxempty$
\begin{remark}[variational approach]
\rm
The key advantage of the variational discretization approach is that delivers an optimal quadratic rate of convergence for the error approximation of the control variable \cite[Theorem 2.4]{MR2122182}. The analysis relies on the following assumption \cite[Assumption 2.3]{MR2122182}:
\[
 \| (\mathbf{S} - \mathbf{S}_{\T}) z \|_{L^2(\Omega)} \lesssim h_{\T}^2 \| z \|_{L^2(\Omega)},
\]
which, in turn, relies on the $H^2(\Omega)$-regularity of the optimal state variable $\bar u$. In our problem, the regularity properties exhibited by $\bar u$ are limited. In fact, Corollary \ref{thm:reg_state_smooth} reveals that $\bar u \in H^{s+1/2-\epsilon}(\Omega)$ for every $\epsilon >0$. As \eqref{eq:example} shows, this is the case even when $\partial \Omega$ is smooth. This reduced regularity feature is responsible for the suboptimal order of convergence in the error estimate \eqref{eq:error_estimate_control_va}.
\end{remark}

\subsection{A fully discrete scheme}
\label{subsec:fd}

In this section, we propose and analyze a fully discrete scheme to approximate the solution of the fractional optimal control problem \eqref{eq:min}--\eqref{eq:control_constraints} by using piecewise constant discretization for the approximation of the control variable and piecewise linear discretization for the approximation of the state variable. To be precise, to discretize the control, we introduce the finite element space of piecewise constant functions over $\T$
\begin{equation}
 \mathbb{W}(\T) = \left\{ v_{\T} \in L^{\infty}( \Omega): {v_{\T}}_{|T} \in \mathbb{P}_0(T) \ \forall T \in \T \right\},
\end{equation}
and the space of discrete admissible controls
\begin{equation}
 \mathbb{Z}_{\mathrm{ad}} (\T) = \mathbb{Z}_{\mathrm{ad}} \cap  \mathbb{W}(\T).
\end{equation}

With this notation at hand, we propose the following fully discrete approximation of the optimal control problem \eqref{eq:min}--\eqref{eq:control_constraints}: Find $\min J (u_{\T},z_{\T})$ subject to the discrete state equation
\begin{equation}
\label{eq:weak_state_discrete2}
\mathcal{A}(u_{\T},v_{\T}) = \langle f + z_{\T}, v_{\T} \rangle \quad \forall v_{\T} \in \V(\T),
\end{equation}
and the control constraints $z_{\T} \in \mathbb{Z}_{\mathrm{ad}}(\T)$, where $J$, $\mathcal{A}$, and $\V(\T)$ are defined as in \eqref{eq:cost_functional}, \eqref{eq:bilinear_form}, and \eqref{eq:defFESpace}, respectively. For notational convenience, we will refer to the previously defined problem as the \emph{fully discrete optimal control problem}.

We define the discrete control-to-state operator $\mathbf{S}_{\T}: \mathbb{Z}_{\mathrm{ad}}(\T) \ni z_{\T} \mapsto u_{\T} \in \mathbb{V}(\T)$,
where $\mathbf{S}_{\T} z_{\T} = u_{\T}$ solves \eqref{eq:weak_state_discrete2}. We also define the optimal adjoint state $\bar p_{\T}$ as the solution to
\begin{equation}
\label{eq:weak_adjoint_discrete}
\mathcal{A}(v_{\T},\bar p_{\T}) = (\bar u_{\T} - u_d , v_{\T} )_{L^2(\Omega)} \quad \forall v_{\T} \in \V(\T).
\end{equation}

We present the following result.
\begin{theorem}[existence, uniqueness and optimality conditions]
The fully discrete optimal control problem has a unique optimal solution $(\bar u_{\T},\bar z_{\T}) \in \V(\T) \times  \mathbb{Z}_{\mathrm{ad}}(\T)$. In addition, the first order optimality condition
\begin{equation}
\label{eq:variational_inequality_discrete}
(\bar p_{\T} + \alpha \bar z_{\T}, z_{\T} - \bar  z_{\T} )_{L^2(\Omega)} \geq 0 \quad \forall  z_{\T} \in  \mathbb{Z}_{\mathrm{ad}}(\T)
\end{equation}
is necessary and sufficient.
\end{theorem}
{\it Proof.}
The proof follows standard arguments \cite{MR2583281}. For brevity, we skip the details.
$\boxempty$

\subsubsection{Auxiliary estimates and variables}
Since it is instrumental in the analysis that we perform, we introduce the $L^2(\Omega)$-orthogonal projection operator \cite[Section 1.6.3]{MR2050138}
\begin{equation}
\Pi_{\T}: L^2(\Omega) \rightarrow \mathbb{W}(\T),
\qquad
(v -\Pi_{\T}v,v_\T)_{L^2(\Omega)} = 0 \quad \forall v_\T \in \mathbb{W}(\T).
\label{eq:orthogonal_projection}
\end{equation}
An important property is that $\Pi_{\T} \mathbb{Z}_{\mathrm{ad}} \subset  \mathbb{Z}_{\mathrm{ad}}(\T)$. In addition, for $1 \leq p \leq \infty$, $\kappa \in (0,1]$, and $v \in W^{\kappa,p}(\Omega)$, we have the error estimate \cite[Proposition 1.135]{MR2050138}
\begin{equation}
 \| v - \Pi_{\T} v\|_{L^p(\Omega)} \lesssim h_{\T}^\kappa |v|_{W^{\kappa,p}(\Omega)}.
 \label{eq:estimate_projection}
\end{equation}

In what follows we introduce two auxiliary variables that are also instrumental to perform an error analysis for the fully discrete optimal control problem. First,
\begin{equation}
\label{eq:q}
q_{\T} \in \V(\T): \quad \mathcal{A}(v_{\T},q_{\T}) = (\bar u  - u_d,v_{\T})_{L^2(\Omega)} \quad \forall v_{\T} \in \V(\T).
\end{equation}
Second,
\begin{equation}
\label{eq:r}
r_{\T} \in \V(\T): \quad \mathcal{A}(v_{\T},r_{\T}) = (u_{\T}( \bar z) - u_d,v_{\T})_{L^2(\Omega)} \quad \forall v_{\T} \in \V(\T),
\end{equation}
where $u_{\T}( \bar z) \in \V(\T)$ solves the discrete problem \eqref{eq:weak_state_discrete} with $z$ replaced by $\bar z$.

\subsubsection{A priori error estimates on smooth domains}
We now derive error estimates for the fully discrete optimal control problem when $\partial \Omega$ is smooth.
\begin{theorem}[error estimate for smooth domains on quasi--uniform meshes]\label{thm:L2controlQuasiUniform}
Let $s \in (0,1)$ and $u_d \in H^{1/2 - s - \epsilon}(\Omega)$, for $\epsilon >0$ arbitrarily small. Let $(\bar u, \bar z)$ and $(\bar u_\T, \bar z_\T)$ be the solutions to the continuous and fully discrete optimal control problems, respectively. Let $\partial \Omega$ be a smooth domain and $\T$ be quasi--uniform. Under the framework of Theorem \ref{thm:reg_control_smooth}, we have the following error estimates: If $s > 1/2$, then
\begin{equation}
\label{eq:error_estimate_control}
 \| \bar z - \bar z_{\T}  \|_{L^2(\Omega)} \lesssim h_{\T}^{1-\epsilon} 
 \left(  \| \bar u  \|_{H^{1/2-s -\epsilon}(\Omega)} + 
 \| u_d  \|_{H^{1/2-s -\epsilon}(\Omega)} +  \| \bar z \|_{H^{1}(\Omega)}  \right)
\end{equation}
and if $s \leq 1/2$, then
\begin{equation}
\label{eq:error_estimate_control_s005}
 \| \bar z - \bar z_{\T}  \|_{L^2(\Omega)} \lesssim h_{\T}^{s + \frac{1}{2} - \epsilon} 
 \left(  \| \bar u  \|_{H^{1/2 -s-\epsilon}(\Omega)} + \| u_d  
 \|_{H^{1/2-s -\epsilon}(\Omega)} + \| \bar z \|_{H^{s + 1/2 - \epsilon}(\Omega)} \right).
\end{equation}
In both estimates the hidden constants depend on $\Omega$, $n$, and $s$.
\label{thm:error_estimate}
\end{theorem}
{\it Proof.}
We proceed in four steps.

\noindent \emph{Step} 1. We begin this step by observing that, since $\mathbb{Z}_{\mathrm{ad}}(\T) \subset \mathbb{Z}_{\mathrm{ad}}$, we are allow to set $z = \bar z_{\T}$ in the optimality condition \eqref{eq:variational_inequality}. On the other hand, we set $z_{\T} = \Pi_{\T} \bar z \in \mathbb{Z}_{\mathrm{ad}}(\T)$ in \eqref{eq:variational_inequality_discrete}; $\Pi_{\T}$ denotes the $L^2(\Omega)$-orthogonal projection operator defined in \eqref{eq:orthogonal_projection}. Adding the obtained inequalities, we arrive at the estimate
\begin{equation}
\label{eq:I+II}
\alpha  \|  \bar z - \bar z_{\T} \|_{L^2(\Omega)} \leq  (\bar p - \bar p_{\T}, \bar z_{\T} - \bar z)_{L^2(\Omega)} +  (\bar p_{\T} + \alpha \bar z_{\T}, \Pi_{\T} \bar z - \bar z)_{L^2(\Omega)} =: \textrm{I} + \textrm{II}.
\end{equation}

\noindent \emph{Step} 2. We bound $\textrm{I}$. To accomplish this task, we write $\bar p - \bar p_{\T} = (\bar p - q_{\T}) + (q_{\T} - \bar p_{\T})$, where $q_{\T}$ is defined as in \eqref{eq:q} and first estimate the term involving $\bar p - q_{\T}$. Since $q_{\T}$ can be seen as the finite element approximation of $\bar p$ within the space $\V(\T)$, we can thus invoke the a priori error estimate \eqref{eq:L2estimater} with $r = 1/2-s-\epsilon$ to conclude the estimate
\begin{equation}
 \| \bar p - q_{\T}  \|_{L^2(\Omega)} \lesssim h_{\T}^{\vartheta_1 + \beta} \left(  \| \bar u  \|_{H^{1/2 -s - \epsilon}(\Omega)} + \| u_d  \|_{H^{1/2 - s - \epsilon}(\Omega)}  \right),
 \label{eq:pminusq}
\end{equation}
where $\vartheta_1 = 1/2 - \epsilon$, $\beta = \min \{ s, 1/2 - \epsilon \}$, and $\epsilon > 0$ being arbitrarily small. Notice that Corollary \ref{thm:reg_state_smooth} guarantees that $\bar u \in H^{s + 1/2 - \epsilon}(\Omega)$. Thus $\bar u \in H^{1/2-s-\epsilon}(\Omega)$ and, by assumption, $\bar u - u_d \in H^{1/2-s-\epsilon}(\Omega)$ for every $\epsilon > 0$.


To control the term $q_{\T} - \bar p_{\T}$, we write $q_{\T} - \bar p_{\T} = (q_{\T} - r_{\T}) + (r_{\T} - \bar p_{\T} )$, where $r_{\T}$ is defined as in \eqref{eq:r}. Next, notice that $r_{\T} - \bar p_{\T} \in \V(\T)$ solves
\begin{equation}
\label{eq:aux1}
\mathcal{A}(v_{\T},r_{\T} - \bar p_{\T}) = (u_{\T}( \bar z) - \bar u_{\T},v_{\T})_{L^2(\Omega)} \quad \forall v_{\T} \in \V(\T).
\end{equation}
On the other hand, $u_{\T}( \bar z) - \bar u_{\T} \in \V(\T)$ solves
\begin{equation}
\label{eq:aux2}
\mathcal{A}(u_{\T}(\bar z) - \bar u_{\T},v_{\T}) = ( \bar z - \bar z_{\T},v_{\T})_{L^2(\Omega)} \quad \forall v_{\T} \in \V(\T).
\end{equation}
Consequently, by setting $v_{\T} = r_{\T} - \bar p_{\T} \in \V(\T)$ in \eqref{eq:aux2} and $v_{\T} = \bar u_{\T} - u_{\T}(\bar z) \in \V(\T)$ in \eqref{eq:aux1}, we conclude that
\begin{equation}
 (r_{\T} - \bar p_{\T}, \bar z_{\T} - \bar z)_{L^2(\Omega)} = \mathcal{A}(\bar u_{\T} - u_{\T}(\bar z), r_{\T} - \bar p_{\T}) = - \| \bar u_{\T} - u_{\T}(\bar z) \|_{L^2(\Omega)}^2 \leq 0.
 \label{eq:negative_sign}
\end{equation}
It thus suffices to estimate $q_{\T} - r_{\T}$. To accomplish this task, we notice that
\begin{equation*}
q_{\T} - r_{\T} \in \V(\T): \quad \mathcal{A}(v_{\T},q_{\T} - r_{\T}) = (\bar u  - u_{\T}(\bar z),v_{\T})_{L^2(\Omega)} \quad \forall v_{\T} \in \V(\T).
\end{equation*}
We invoke a stability argument and the a priori error estimate \eqref{eq:L2estimater} with $r = \gamma = \min \{ s+1/2-\epsilon,1\}$ to conclude that
\begin{equation}
 \| q_{\T} - r_{\T}  \|_{L^2(\Omega)} \lesssim  \| \bar u  - u_{\T}(\bar z) \|_{L^2(\Omega)} \lesssim h_{\T}^{\vartheta_2 + \beta}  \| \bar z \|_{H^{\gamma}(\Omega)},
\label{eq:qminusr}
\end{equation}
where $\vartheta_2 = \min \{ s + \gamma, 1/2 - \epsilon\}$, and $\beta = \min \{ s,1/2 - \epsilon\}$. The fact that $\bar z \in H^{\gamma}(\Omega)$ follows from Theorem \ref{thm:reg_control_smooth}.

In view of Young's inequality, the collection of the estimates  \eqref{eq:pminusq}, \eqref{eq:negative_sign}, and \eqref{eq:qminusr} yield the estimate for the term $\mathrm{I}$:
\begin{multline*}
|\mathrm{I}|
 \leq C h_{\T}^{2( \vartheta_1 + \beta)} \left(  \| \bar u  \|^2_{H^{1/2-s-\epsilon}(\Omega)} + \| u_d  \|^2_{H^{1/2-s-\epsilon}(\Omega)} \right)
 \\
+  C h_{\T}^{2( \vartheta_2 + \beta)}  \| \bar z \|^2_{H^{\gamma}(\Omega)}
 + \frac{\alpha}{4}  \| \bar z  - \bar z_{\T}\|_{L^2(\Omega)}^2,
\end{multline*}
where $C$ denotes a positive constant. We note that, for $s \in (0,1)$, $\vartheta_1 = 1/2 - \epsilon$ and $\vartheta_2 = 1/2 - \epsilon$. We can thus conclude the following estimates for the term $\mathrm{I}$:
\begin{align}
s \in [1/2,1)  \Rightarrow \beta = 1/2 - \epsilon  & \Rightarrow |\mathrm{I}|  \lesssim C h_{\T}^{2(1-\epsilon)},
\\
s \in (0,1/2)  \Rightarrow \beta = s  & \Rightarrow |\mathrm{I}|  \lesssim C h_{\T}^{2(s+1/2-\epsilon)}.
\end{align}

\noindent \emph{Step} 3.  The goal of this step is to estimate the term $\mathrm{II} = (\bar p_{\T} + \alpha \bar z_{\T}, \Pi_{\T} \bar z - \bar z)_{L^2(\Omega)}$. To accomplish this task, we invoke definitions \eqref{eq:q} and \eqref{eq:r} and write
\begin{multline*}
\mathrm{II} = (\bar p_{\T} + \alpha \bar z_{\T}, \Pi_{\T} \bar z - \bar z)_{L^2(\Omega)} = (\bar p + \alpha \bar z, \Pi_{\T} \bar z - \bar z)_{L^2(\Omega)}
+ \alpha (\bar z_{\T} - \bar z, \Pi_{\T} \bar z - \bar z)_{L^2(\Omega)}
\\
+  (\bar p_{\T} - r_{\T}, \Pi_{\T} \bar z - \bar z)_{L^2(\Omega)} +  (r_{\T} \pm q_{\T} - \bar p, \Pi_{\T} \bar z - \bar z)_{L^2(\Omega)} =: \mathrm{II}_1 +   \mathrm{II}_2 + \mathrm{II}_3 + \mathrm{II}_4.
\end{multline*}
To bound $\mathrm{II}_1$, we first invoke the definition of $\Pi_{\T}$ and notice that
\begin{equation}
 \mathrm{II}_1 = (\bar p + \alpha \bar z - \Pi_{\T}(\bar p + \alpha \bar z), \Pi_{\T} \bar z - \bar z)_{L^2(\Omega)}.
 \label{eq:II_1}
\end{equation}
We can thus invoke the estimate \eqref{eq:estimate_projection} and the regularity results of Theorem \ref{thm:reg_control_smooth} and Corollary \ref{thm:reg_state_smooth} to conclude that
 \begin{equation}
| \mathrm{II}_1| \lesssim h_{\T}^{2\gamma} \| \bar p  + \alpha  \bar z \|_{H^{\gamma}(\Omega)} \| \bar z \|_{H^{\gamma}(\Omega)},
 \end{equation}
where $\gamma = \min \{ s+1/2 -\epsilon, 1 \}$ with $\epsilon >0$ arbitrarily small.

To bound $ \mathrm{II}_2$ we use, \eqref{eq:estimate_projection} and the regularity results of Theorem \ref{thm:reg_control_smooth}, again, and Young's inequality. We thus arrive at the estimate
  \begin{equation}
  | \mathrm{II}_2| \leq C h_{\T}^{2\gamma} \| \bar z \|^2_{H^{\gamma}(\Omega)} + \frac{\alpha}{4} \| \bar z  -  \bar z_{\T} \|^2_{L^2(\Omega)},
 \end{equation}
where $C$ denotes a positive constant and $\gamma = \min \{ s+1/2 -\epsilon, 1 \}$.

To control $ \mathrm{II}_3$ we invoke a stability estimate for the discrete problem \eqref{eq:r} and the error estimate \eqref{eq:estimate_projection}. In fact, we have that
\begin{equation}
  | \mathrm{II}_3| \lesssim h_{\T}^{\gamma} \|  \bar u_{\T} - u_{\T}(\bar z) \|_{L^{2}(\Omega)} \|  \bar z \|_{H^{\gamma}(\Omega)}
  \lesssim h_{\T}^{\gamma} \| \bar z - \bar{z}_{\T} \|_{L^{2}(\Omega)} \|  \bar z \|_{H^{\gamma}(\Omega)} ,
\end{equation}
To obtain the last inequality we have used a stability estimate for the discrete problem \eqref{eq:weak_state_discrete2}.

The control of the term $r_{\T} - q_{\T}$ follows from \eqref{eq:qminusr} while the one for $q_{\T} - \bar p$ from the estimate \eqref{eq:pminusq}.

\noindent \emph{Step} 4. The desired estimates \eqref{eq:error_estimate_control} and \eqref{eq:error_estimate_control_s005}
follow from collecting all the estimates we obtained in previous steps.
$\boxempty$

\subsubsection{A priori error estimates on Lipschitz domains}
To derive the error estimates \eqref{eq:error_estimate_control} and \eqref{eq:error_estimate_control_s005} we have used the a priori error estimate \eqref{eq:L2estimater} that requires that $\partial\Omega$ is smooth. In the following result we allow $\Omega$ to be a bounded Lipschitz domain satisfying the exterior ball condition and obtain quasi-optimal error estimates, in terms of approximation, for the control and state variables. To do this, we define
\begin{equation}
 \Lambda(\bar z, f, u_d):= \| f + \bar z \|_{C^{1-s}(\overline \Omega)}
 +  \| u_d \|_{C^{1-s}(\overline \Omega)} +  \| \bar z \|_{H^1(\Omega)}.
\end{equation}
We present the following result.
\begin{theorem}[error estimates for Lipschitz domains on graded meshes]\label{thm:errorGraded}
Let $s \in (1/2,1)$ and $\Omega$ be a bounded Lipschitz domain satisfying the exterior ball condition. Let $(\bar u, \bar z)$ and $(\bar u_\T, \bar z_\T)$ be the solutions to the continuous and fully discrete optimal control problems, respectively. If $\T$ satisfies \eqref{eq:graded_meshes} with $\mu = 2$, $f \in C^{1-s}(\overline\Omega)$, and $u_d \in C^{1-s}(\overline \Omega)$, then
\begin{equation}
\label{eq:error_estimate_control2}
 \| \bar z - \bar z_{\T}  \|_{L^2(\Omega)} \lesssim |\log N| N^{-\frac{1}{2}}  \Lambda(\bar z, f, u_d)
\end{equation}
and
\begin{equation}
\label{eq:error_estimate_state}
 \| \bar u - \bar u_{\T}  \|_{s} \lesssim |\log N| N^{-\frac{1}{2}} \Lambda(\bar z, f, u_d),
\end{equation}
\label{thm:error_estimate2}
where $N$ denotes the number of degrees of freeedom of $\T$. In both estimates, the hidden constant depend on $\sigma$ and blows up when $s \rightarrow 1/2$.
\end{theorem}
{\it Proof.}
The proof follows closely the arguments developed in the proof of Theorem \ref{thm:error_estimate}; the difference being the use of the error estimate \eqref{eq:energy_estimate_graded} instead of \eqref{eq:L2estimater}. Since the latter estimates require different assumptions on the problem data, we briefly report the arguments.

\noindent \emph{Step} 1. We recall the estimate \eqref{eq:I+II}:
\begin{equation}
\label{eq:I+II2}
\alpha  \|  \bar z - \bar z_{\T} \|^2_{L^2(\Omega)} \leq  (\bar p - \bar p_{\T}, \bar z_{\T} - \bar z)_{L^2(\Omega)} +  (\bar p_{\T} + \alpha \bar z_{\T}, \Pi_{\T} \bar z - \bar z)_{L^2(\Omega)} =: \textrm{I} + \textrm{II}.
\end{equation}

\noindent \emph{Step} 2. The results of Proposition \ref{pro:state_regularity_1} imply that $\bar u \in C^{s}(\overline \Omega)$ with the stability estimate $\| \bar u \|_{C^s(\overline \Omega)} \lesssim \| f + \bar z\|_{L^{\infty}(\Omega)}$. This, in view of the assumption $u_d \in C^{1-s}(\overline \Omega)$, allows us to conclude that $\bar u - u_d \in C^{1-s}(\overline \Omega)$ for $s \in (1/2,1)$. Notice that $\bar u - u_d$ corresponds to the right-hand side of the adjoint equation \eqref{eq:weak_adjoint_equation} and that $q_{\T}$ denotes its finite element approximation. We can thus conclude, on the basis of the error estimate \eqref{eq:energy_estimate_graded}, that
\begin{align}
 \| \bar p - q_{\T}  \|_{L^2(\Omega)}
 & \lesssim |\log N| N^{-\frac{1}{2}} \| \bar u - u_d \|_{C^{1-s}(\overline \Omega)}
 \nonumber
 \\
 & \lesssim |\log N| N^{-\frac{1}{2}} \left( \| f + \bar z\|_{L^{\infty}(\Omega)} + \|  u_d \|_{C^{1-s}(\overline\Omega)} \right).
 \label{eq:pminusq2}
\end{align}

The control of $q_{\T} - \bar p_{\T}$ follows from writing $q_{\T} - \bar p_{\T} = (q_{\T} - r_{\T}) + (r_{\T} - \bar p_{\T} )$, where $r_{\T}$ is defined as in \eqref{eq:r}. Notice that \eqref{eq:negative_sign} yields
\begin{equation}
 (r_{\T} - \bar p_{\T}, \bar z_{\T} - \bar z)_{L^2(\Omega)} = - \| \bar u_{\T} - u_{\T}(\bar z) \|_{L^2(\Omega)}^2 \leq 0.
 \label{eq:negative_sign2}
\end{equation}
Now, notice that, in view of \eqref{eq:regularity_z_as_rhs}, the optimal control $\bar z \in C^{s}(\overline \Omega)$ when $s \in (1/2,1)$. Consequently, for such an interval, $f + \bar z \in C^{1-s}(\overline\Omega)$. We thus invoke a stability argument and the error estimate \eqref{eq:energy_estimate_graded} to conclude that
\begin{equation}
 \| q_{\T} - r_{\T}  \|_{L^2(\Omega)} \lesssim  \| \bar u  - u_{\T}(\bar z) \|_{s} \lesssim |\log N|  N^{-\frac{1}{2}} \| f + \bar z \|_{C^{1-s}(\overline\Omega)}.
\label{eq:qminusr2}
\end{equation}

\noindent \emph{Step} 3. As in the step 3 in the proof of Theorem \ref{thm:error_estimate}, we write $\textrm{II} = \textrm{II}_1 + \textrm{II}_2 + \textrm{II}_3 + \textrm{II}_4$.  The estimate for $\textrm{II}_1$ follows from \eqref{eq:II_1} and the error estimate \eqref{eq:estimate_projection}:
\begin{equation}
 |\mathrm{II}_1| \lesssim h^2 \| \bar p + \alpha \bar z \|_{H^1(\Omega)} \| \bar z \|_{H^1(\Omega)},
\end{equation}
where we have used that the mesh grading \eqref{eq:graded_meshes} implies that $h_T \leq C h$ for all $T \in \T$. Notice that, in view of the regularity estimates of Theorem \ref{thm:reg_control_2} we have that $\| \bar p \|_{H^1(\Omega)} $ and $\| \bar z \|_{H^1(\Omega)}$ are bounded. The estimate for the term $\mathrm{II}_2$ follows from the regularity estimates of Theorem \ref{thm:reg_control_2} and the error estimate \eqref{eq:estimate_projection}:
  \begin{equation}
  | \mathrm{II}_2| \leq C h^{2} \| \bar z \|^2_{H^{1}(\Omega)} + \frac{\alpha}{4} \| \bar z  -  \bar z_{\T} \|^2_{L^2(\Omega)},
 \end{equation}
where $C$ denotes a positive constant. The estimates for $\textrm{II}_3$ and $\textrm{II}_4$ follow form the estimates derived for $r_{\T}  -  q_{\T}$ and $q_{\T}  -  p$.

\noindent \emph{Step} 4. The desired estimate \eqref{eq:error_estimate_control2} follows from collecting the estimates derived in the previous steps.

\noindent \emph{Step} 5. We derive the error estimates associated to the approximation of the optimal state variable. We begin with the basic estimate
\begin{align*}
 \|  \bar u - \bar u_{\T} \|_{s} & =  \|  \mathbf{S} \bar z - \mathbf{S}_{\T} \bar z_{\T} \|_{s}
 \\
 & \leq
 \|  (\mathbf{S} -  \mathbf{S}_{\T} ) \bar z   \|_{s} +  \|  \mathbf{S}_{\T} (\bar z  -  \bar z_{\T}) \|_{s}.
\end{align*}
Notice that \eqref{eq:regularity_z_as_rhs} guarantees that $\bar z \in C^{s}(\overline\Omega)$ for $s \in (\tfrac{1}{2},1)$ and thus that $\bar z \in C^{1-s}(\overline\Omega)$. We can thus apply the error estimate \eqref{eq:energy_estimate_graded} to conclude that
\[
 \|  (\mathbf{S} -  \mathbf{S}_{\T} ) \bar z   \|_{s} \lesssim |\log N| N^{-\frac{1}{2}} \|  f + \bar z \|_{C^{1-s}(\overline\Omega)},
\]
where the hidden constant depends on $\sigma$ and blows up when $s \rightarrow 1/2$. We now invoke the continuity of the discrete control-to-state map $\mathbf{S}_{\T}$ to conclude that
\[
\|  \mathbf{S}_{\T} (\bar z  -  \bar z_{\T}) \|_{s} \lesssim |\log N| N^{-\frac{1}{2}} \Lambda(\bar z, f, u_d).
\]
The collection of these estimates yield \eqref{thm:error_estimate2}.
$\boxempty$

\begin{remark}[quasi-optimal error estimate]\rm
Notice that the error estimates \eqref{eq:error_estimate_control2} and \eqref{eq:error_estimate_state} are quasi-optimal in terms of approximation.
\end{remark}

\section{Numerical experiments}\label{sec:num_experiments}
We present a series of numerical examples that illustrate the performance of the fully discrete scheme proposed in section \ref{subsec:fd} for the solution of the optimal control problem \eqref{eq:min}--\eqref{eq:control_constraints} and the sharpness of the derived error estimates. We consider an example where $\Omega$ is smooth and another one where we go beyond the theory and violate
the assumption of exterior ball condition.

When solving equations involving the integral fractional Laplacian, two primary issues need to be addressed:
\begin{itemize}
\item
  No closed form is available for the entries of the stiffness matrix, and hence quadrature needs to be used for their evaluation.
  Particular care in the choice of quadrature rules needs to be taken to handle the case of pairs of elements that are either connected or close to each other.
  In order not to spoil the solution, the quadrature error needs to be smaller than the error arising from discretization.
\item
  Due to the nonlocal interactions, straightforward assembly would lead to a dense matrix representation of the fractional Laplacian.
  This would mean that a single solve of state or adjoint equation would scale at best quadratically in the number of unknowns.
  Fortunately, the interactions of well-separated clusters of unknowns can be approximated, using a panel clustering approach, whereby the overall complexity of a matrix-vector product is reduced to \(\mathcal{O}(N\left(\log N\right)^{2n})\).
  Again, error due to the approximation of the operator needs to be controlled.
\end{itemize}
For a comprehensive treatment of both issues we refer the reader to
\cite{AinsworthGlusa2017_AspectsAdaptiveFiniteElement,AinsworthGlusa2018_TowardsEfficientFiniteElement}.

For the examples that we present in this section, the discrete equations \eqref{eq:weak_state_discrete2} and \eqref{eq:weak_adjoint_discrete} are solved on the basis of multigrid solver, while to solve the minimization problem, we use the BFGS algorithm~\cite{NocedalWright2006_NumericalOptimization2nd}.

\subsection{Unit disc}
We let $n=2$, $\Omega = B(0,1)$, and $s \in (0,1)$. We consider
\begin{equation*}
(-\Delta)^s u = f \; \text{ in }\Omega, \qquad u=0 \; \text{ in }\Omega^{c}.
\end{equation*}
This problem has a family of known closed-form solutions when the right-hand side reads, in polar coordinates, as follows:
\begin{equation*}
  f_{n,\ell}(r,\theta)
  = 2^{2s}\Gamma\left(1+s\right)^{2}\binom{s+n+\ell}{s}\binom{s+n}{s} r^{\ell}\cos\left(\ell\theta\right) P_{n}^{(s,\ell)}\left(2r^{2}-1\right),
\end{equation*}
where $\ell,n \in \mathbb{N}_0$. In fact, for $\ell,n \in \mathbb{N}_0$, the solution is given by
\begin{align*}
u_{n,\ell} (r,\theta) = r^{\ell}\cos(\ell\theta) P_{n}^{(s,\ell)}(2r^{2}-1)(1-r^{2})^{s}_{+}.
\end{align*}
We refer the reader to \cite{DydaKuznetsovEtAl2016_FractionalLaplaceOperatorMeijerGFunction} for details.

We set $a = -0.9$, $b = 0.9$, \(\alpha=10^{-1}\), \(u_{d}=u_{0,1}+\alpha f_{0,0}\), and \(f=f_{0,1}-\textrm{proj}_{[a,b]}\left(u_{0,0}\right)\). The exact solution reads \(\bar u=u_{0,1}\), \(\bar p=-\alpha u_{0,0}\) and
\begin{align*}
\bar z&=\textrm{proj}_{[a,b]}\left(u_{0,0}\right) =
     \begin{cases}
     b & r<r_{o}:=\sqrt{1-b^{1/s}}, \\
     (1-r^{2})^{s} & r\geq r_{o}.
     \end{cases}
\end{align*}
\subsubsection{Quasi-uniform meshes}
We discretize $\Omega$ using a sequence of quasi-uniform meshes and solve the control problem with the scheme of section \ref{subsec:fd} for \(s\in\{0.1,0.2,\dots,0.9\}\).
\begin{figure}
\centering
\includegraphics{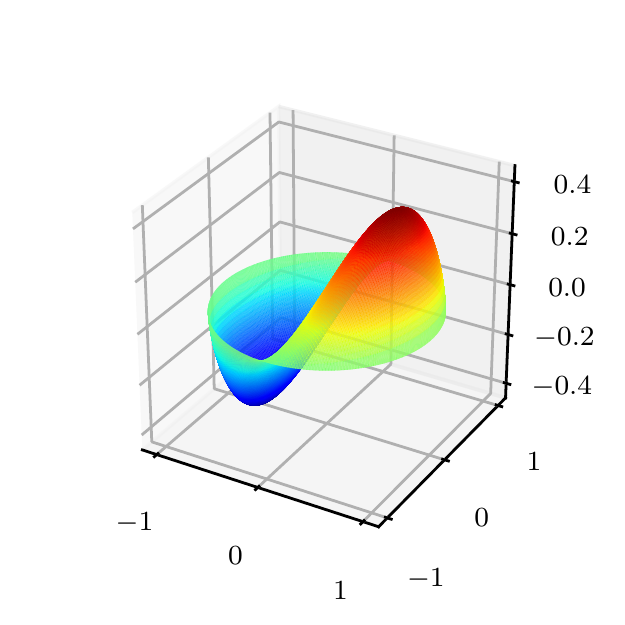}
\includegraphics{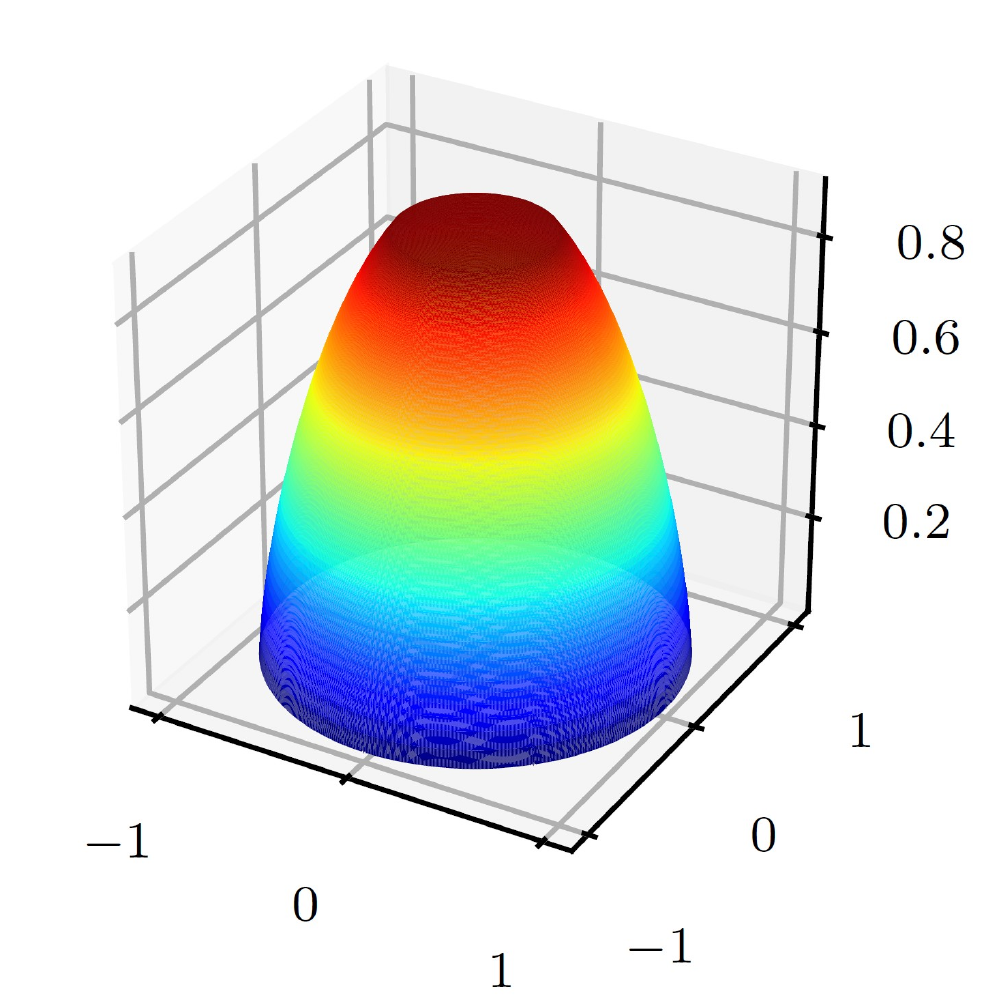}
\caption{Finite element solutions for the optimal state \( \bar u_{\T}\) (\emph{left}) and the optimal control \( \bar z_{\T}\) (\emph{right}) for \(s=0.7\). We notice that the upper bound on the control is active near the center of the domain.}
\label{fig:solutions0.7}
\vspace{-0.5cm}
\end{figure}
In Figure~\ref{fig:solutions0.7} we present the finite element solutions for the optimal state \( \bar u_{\T}\) and control \(\bar z_{\T}\), on the finest mesh (66k vertices, 131k elements), for \(s=0.7\).
Note that the upper bound on the control is active for \(r\leq r_{o}\).

\begin{figure}
\centering
\includegraphics{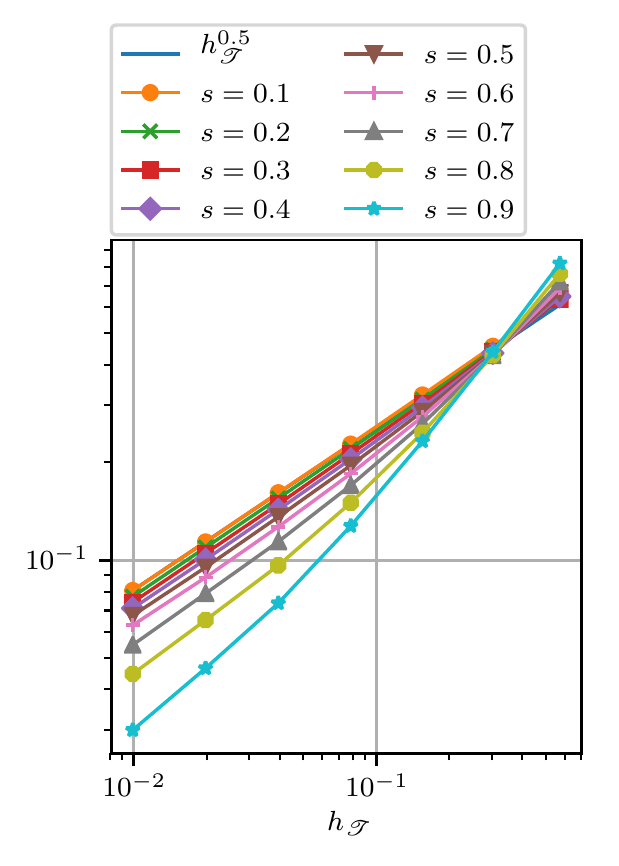}
\includegraphics{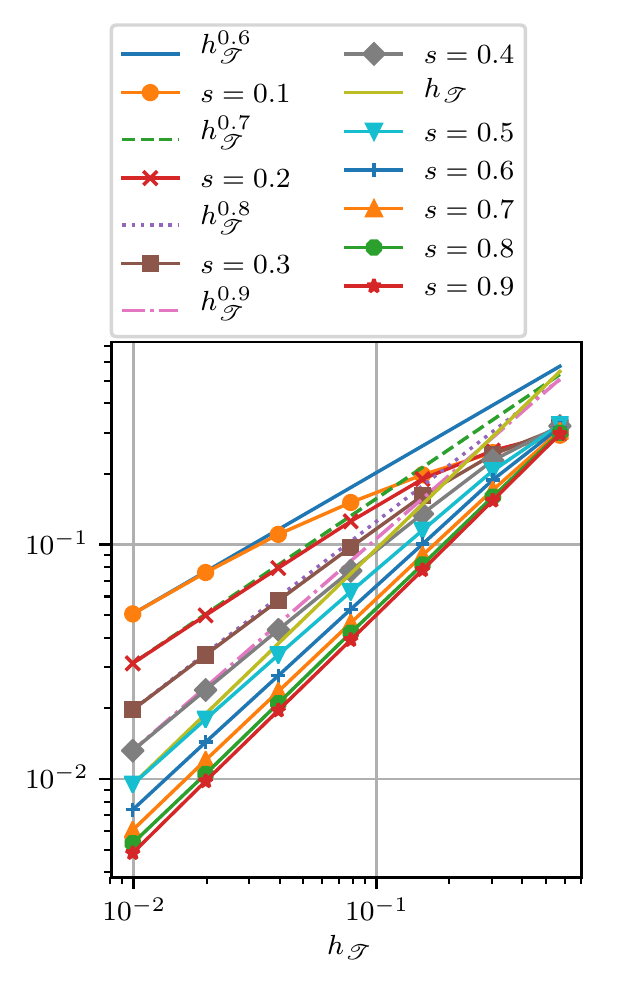}
\caption{\emph{Left:} Experimental rates of convergence for the
\(\tilde{H}^{s}(\Omega)\)-error of the state variable and the \(L^{2}(\Omega)\)-error of the control variable for $n=2$, $\Omega = B(0,1)$, and $s \in \{ 0.1,0.2,\dots,0.9\}$. The experimental rates of convergence are in agreement with the results of  Proposition~\ref{prop:HsstateQuasiUniform} and Theorem~\ref{thm:L2controlQuasiUniform}.}
\label{fig:errors}
\end{figure}

In Figures~\ref{fig:errors}, we show experimental rates of convergence for the \(\tilde{H}^{s}(\Omega)\)-error of the state variable, as well as the \(L^{2}(\Omega)\)-error of the control variable. We mention that the aforementioned \(\tilde{H}^{s}(\Omega)\)-error can be computed as follows:
\begin{align}
\|\bar u-\bar u_{\T}\|_s^{2}&= \mathcal{A}(\bar u- \bar u_{\T},\bar u-\bar u_{\T}) \nonumber
                             = \mathcal{A}(\bar u,\bar u) -2 \mathcal{A}(\bar u,\bar u_{\T})+ 
                               \mathcal{A}(\bar u_{\T},\bar u_{\T}) \nonumber \\
                            &= \langle f + \bar z,\bar u \rangle - 2\langle f + 
                               \bar z,\bar u_{\T} \rangle + \langle f + 
                               z_{\T},\bar u_{\T} \rangle \nonumber \\
                            &= \langle f + \bar z,\bar u \rangle - 2\langle \bar z, 
                               \bar u_{\T} \rangle - \langle f,\bar u_{\T} \rangle 
                               + \langle z_{\T},\bar u_{\T} \rangle, \label{eq:numericalHs}
\end{align}
where the first term can be evaluated analytically. We observe, from Figures~\ref{fig:errors}, that the rates of convergence predicted by Proposition~\ref{prop:HsstateQuasiUniform} and Theorem~\ref{thm:L2controlQuasiUniform} are attained: we observe
\( \mathcal{O}(h_{\T}^{1/2-\epsilon}) \) for the \(\tilde{H}^{s}(\Omega)\)-error of the state variable, and
\[
 \mathcal{O}(h_{\T}^{s+1/2-\epsilon}) \textrm{ and } \mathcal{O}(h_{\T}^{1-\epsilon}),
\]
for the \(L^{2}(\Omega)\)-error of control variable when $s\leq 1/2$ and $s>1/2$, respectively.

\begin{figure}
\centering
\includegraphics{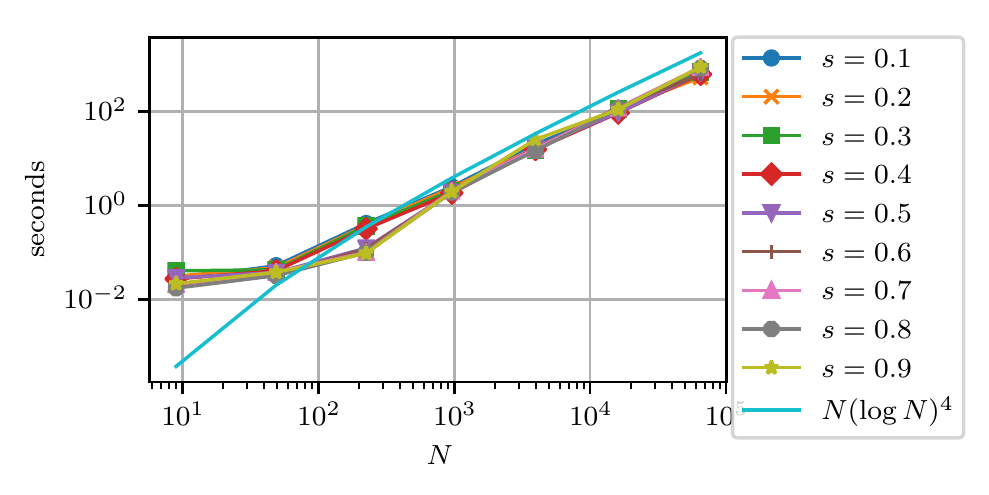}
\caption{Solution time for the discretized control problem using BFGS for the minimization problem and multigrid combined with panel clustering for the linear systems. The solve scale as \(N(\log N)^{4}\).}
\label{fig:timings}
\end{figure}

Figure~\ref{fig:timings} displays the solution times for the discretized control problems. It can be observed that the solve in fact scales as \(\mathcal{O}\left(N(\log N)^{4}\right)\).

\subsubsection{Graded meshes}
\begin{figure}
\centering
\includegraphics[scale=0.8]{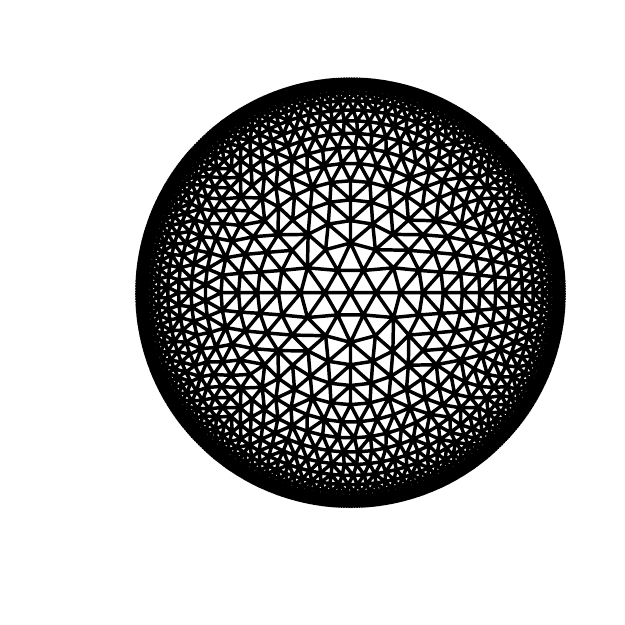}
\includegraphics{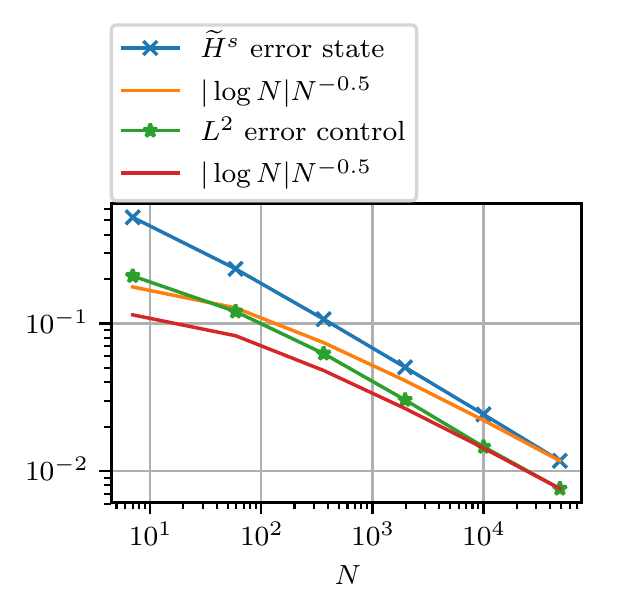}
\caption{\emph{Left:} Graded mesh satisfying condition~\eqref{eq:graded_meshes} with \(\mu=2\). \emph{Right:} Experimental rates of convergence for the \(\tilde{H}^{s}(\Omega)\)-error for the state variable and the \(L^{2}(\Omega)\)-error for the control variable. As predicted by Theorem~\ref{thm:errorGraded}, both experimental rates decay as \(\mathcal{O}(|\log N| N^{-\frac{1}{2}})\), which is nearly optimal in terms of approximation.}
\vspace{-0.5cm}
\label{fig:graded}
\end{figure}

We discretize $\Omega$ using a family of graded meshes which satisfy condition~\eqref{eq:graded_meshes} with \(\mu=2\). As an example, we present one of these meshes in Figure~\ref{fig:graded}. We solve the fractional optimal control problem for \(s=0.75\).
In Figure~\ref{fig:graded}, we present the experimental orders of convergence for the \(\tilde{H}^{s}(\Omega)\)-error for the state variable and the \(L^{2}(\Omega)\)-error for the control variable; both of them being displayed versus the number of degrees of freedom $N$, where, we recall that, $N=\dim\V(\T)$.
It can be observed that, as predicted by Theorem~\ref{thm:errorGraded}, the experimental errors decay as \(\mathcal{O}(|\log N| N^{-\frac{1}{2}})\); the latter being nearly-optimal in terms of approximation.

\FloatBarrier
\subsection{L-shaped domain}
We now illustrate the case of a non-smooth domain by solving the fractional optimal control problem on a family of quasi-uniform meshes on the L-shaped domain \(\Omega=[0,2]^{2}\setminus [1,2]^{2}\). Notice that $\Omega$ is Lipschitz but does not satisfy the exterior ball condition.

\begin{figure}
\centering
\includegraphics{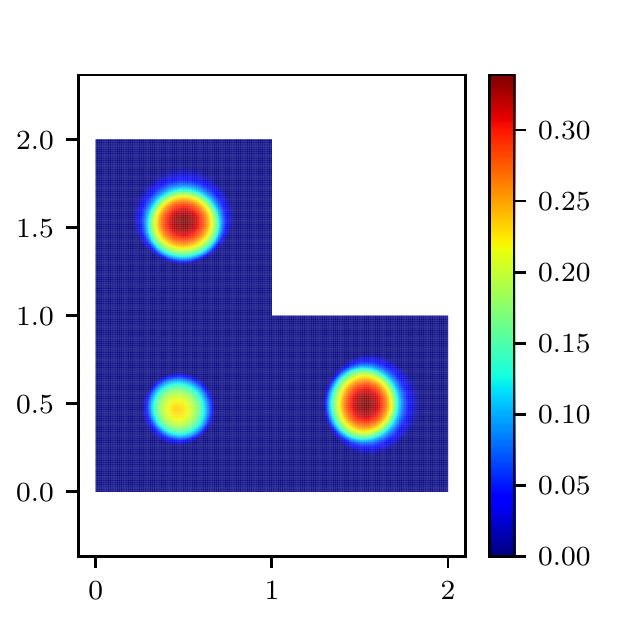}
\includegraphics{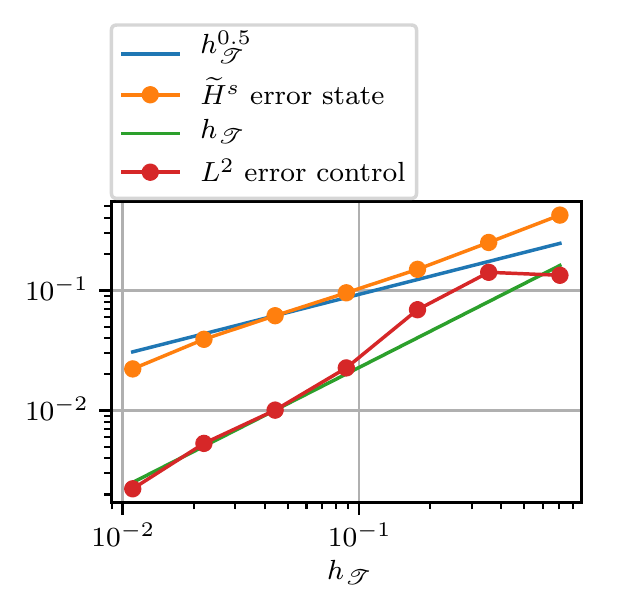}
\caption{\emph{Left:} Finite element solution for the optimal control \( \bar z_{\T}\). \emph{Right:} Experimental rates of convergence for the \(\tilde{H}^{s}(\Omega)\)-error of the state variable and the \(L^{2}(\Omega)\)-error of  the control variable for $n=2$, $s =0.75$, and $\Omega = [0,2]^{2}\setminus [1,2]^{2}$. The experimental convergence rates are in agreement with the results of Proposition~\ref{prop:HsstateQuasiUniform} and Theorem~\ref{thm:L2controlQuasiUniform}: \(\mathcal{O}(h_{\T}^{1/2-\epsilon})\) and \(\mathcal{O}(h_{\T}^{1-\epsilon})\), respectively.}
\label{fig:Lshape}
\end{figure}

We consider $s = 0.75$, \(u_{d}=\mathbf{1}_{B((0.5,0.5),0.2)}+\mathbf{1}_{B((1.5,0.5),0.2)}+\mathbf{1}_{B((0.5,1.5),0.2)}\), \(f=1\), $a= 0$, $b=30$, and \(\alpha=10^{-1}\).
Since no analytical solution is available, we compute errors with respect to a reference solution on a highly refined mesh (200k vertices, 400k elements, \(h=2^{-8}\)).
The numerical solution for the control as well as computed errors are shown in Figure~\ref{fig:Lshape}.
The speed-up of convergence in \(\tilde{H}^{s}(\Omega)\)- and \(L^{2}(\Omega)\)-norm for larger number of unknowns is due to the fact that the reference solution is used in their computation instead of the true solution.

\section{Conclusion}
In this paper we introduced an optimal control problem for the integral form of the fractional Laplacian operator with the goal of determining the optimal source term such that the nonlocal solution is as close as possible to a given data. We performed a careful and detailed mathematical and numerical analysis proving well-posedness of the control problem and establishing resularity estimates and convergence results for two finite-dimensional approximations of the continuous problem. Also, we provided several two-dimensional numerical results that illustrate the theory and additional results on complex geometries that show applicability of our approach to more realistic problems.

This work sets the ground for future research: as an example, one could consider a different control variable such as a diffusion parameter or the fractional order itself. The latter problem is very challenging both in terms of analysis (for different controls the solution belongs to a different functional space) and computations (the matrix of the discretized problem needs to be reassembled at each iteration of the optimization algorithm).

\section{Acknowledgments}
Enrique Ot\'arola was supported by CONICYT through FONDECYT project 3160201. Marta D'Elia and Christian Glusa were supported by Sandia National Laboratories (SNL), SNL is a multimission laboratory managed and operated by National Technology and Engineering Solutions of Sandia, LLC., a wholly owned subsidiary of Honeywell International, Inc., for the U.S. Department of Energy’s National Nuclear Security Administration contract number DE-NA-0003525. This paper describes objective technical results and analysis. Any subjective views or opinions that might be expressed in the paper do not necessarily represent the views of the U.S. Department of Energy or the United States Government. SAND Number: SAND2018-11499 O.

\end{document}